\documentclass[ijds,sglanonrev]{informs4}
\usepackage{eqndefns-left} 
\RequirePackage{tgtermes}
\RequirePackage{newtxtext}
\RequirePackage{newtxmath}
\RequirePackage{bm}
\RequirePackage{endnotes}

\OneAndAHalfSpacedXII 

\usepackage{algorithm}
\usepackage{algpseudocode}
\usepackage{tikz}
\usepackage{amsmath}
\usepackage{xr}
\usepackage[hidelinks]{hyperref}
\makeatletter
\newcommand*{\addFileDependency}[1]{
\typeout{(#1)}
\IfFileExists{#1}{}{\typeout{No file #1.}}
}\makeatother

\allowdisplaybreaks

\usepackage{natbib}
 \bibpunct[, ]{(}{)}{,}{a}{}{,}%
 \def\bibfont{\small}%

\EquationsNumberedThrough   
\TheoremsNumberedThrough  
\ECRepeatTheorems  %

\begin{document}

\RUNAUTHOR{Stratman, Boutilier, and Albert}

\RUNTITLE{Decision-Aware Predictive Model Selection for Workforce Allocation}

\TITLE{Decision-Aware Predictive Model Selection for Workforce Allocation}

\ARTICLEAUTHORS{
\AUTHOR{Eric G. Stratman}
\AFF{Department of Industrial and Systems Engineering,
University of Wisconsin-Madison, \EMAIL{egstratman@wisc.edu}}

\AUTHOR{Justin J. Boutilier}
\AFF{Telfer School of Management,
University of Ottawa, \EMAIL{boutilier@telfer.uOttawa.ca}}

\AUTHOR{Laura A. Albert}
\AFF{Department of Industrial and Systems Engineering,
University of Wisconsin-Madison, \EMAIL{laura@engr.wisc.edu}}
}

\ABSTRACT{
Many organizations depend on human decision-makers to make subjective decisions, especially in settings where information is scarce. Although workers are often viewed as interchangeable, the specific individual assigned to a task can significantly impact outcomes due to their unique decision-making processes and risk tolerance. In this paper, we introduce a novel framework that utilizes machine learning to predict worker behavior and employs integer optimization to strategically assign workers to tasks. Unlike traditional methods that treat machine learning predictions as static inputs for optimization, in our approach, the optimal predictive model used to represent a worker's behavior is determined by how that worker is allocated within the optimization process. We present a decision-aware optimization framework that integrates predictive model selection with worker allocation. Collaborating with an auto-insurance provider and using real-world data, we evaluate the effectiveness of our proposed method by applying three different techniques to predict worker behavior. Our findings show the proposed decision-aware framework outperforms traditional methods and offers context-sensitive and data-responsive strategies for workforce management.
}


\maketitle

\section{Introduction}
Humans supplement known information with intuition and experience to make decisions. This ability enables them to make informed decisions, even in situations with limited data and uncertain outcomes. As a result, despite the rise of artificial intelligence and automation, many organizations still rely on humans to assess value and risk in decision-making. For example, a mortgage underwriter at a bank decides whether to approve or deny a loan application; a doctor in a medical practice advises a patient on undergoing a procedure; and a hiring manager in a corporation determines whether to hire a candidate. These decisions inherently rely on a worker's subjective judgment, which complements objective information. 

Although uniform decision-making within a workforce might seem ideal, a worker's personal judgment and experience invariably shape their decision-making process. For example, one poor decision may alter a worker's behavior when faced with a similar decision in the future. Similarly, psychological studies have shown that personality influences risk aversion, making outcomes highly dependent on the worker assigned to a particular decision \citep{becker_relationship_2012, lonnqvist_measuring_2015,el_othman_personality_2020}. Recognizing this, aligning task assignments with a worker's unique characteristics and experience can enhance both the effectiveness and efficiency of an organization's decision-making processes.

In this paper, we explore worker individuality in subjective decisions and investigate leveraging worker tendencies in a workforce allocation problem. We aim to uncover how unique worker tendencies can be effectively utilized in the assignment of tasks, focusing on aligning tasks with workers most likely to make beneficial decisions. Existing research that studies worker allocation typically assumes that worker behavior is uniform or that worker attributes are known with certainty \citep{de_bruecker_workforce_2015}. However, in many real-world applications, behavior varies between individuals and must be learned from data. We examine the challenge of learning worker behavior using real-world data and explore how different approaches impact the resulting workforce allocation problem. 

The framework we present leverages machine learning to model worker behavior in decision-making tasks in combination with integer optimization to assign workers to tasks. The goal is to improve decision-making outcomes by aligning task assignments with each worker's unique characteristics, while supporting the organization's broader objectives. Our framework is not intended to replace human decision-makers. Instead, our approach aims to complement human judgment by aligning task assignments with workers’ unique traits, thus supporting the organization’s broader objectives without removing individual decision-making nuances. For example, if an organization wishes to reduce risk, the optimization model may assign a high-risk task to a worker who has demonstrated risk aversion to similar tasks in the past. Worker assignments are therefore highly dependent on the insights gained from the predictive model. 

The predictive models used to represent worker behavior and their connection to the optimization task are an important consideration. Several approaches exist for incorporating predictive models into optimization tasks, each varying in complexity and effectiveness. The simpliest approach -- the Decoupled Prediction Optimization (DPO) method -- treats prediction and optimization separately. Predictions are made as single-point estimates and then used as deterministic parameters within the optimization model. While DPO's simplicity is appealing, the major drawback is that it does not consider the accuracy of the predictions, potentially leading to less robust decisions.

Uncertainty-Aware Optimization (UAO) takes a more nuanced approach by accounting for prediction error or parameter uncertainty within the optimization model. UAO typically employs robust or stochastic optimization methods to account for uncertainty. This approach aims to generate solutions that are resilient to the uncertainty inherent in predictive models, thus enhancing decision quality. However, uncertainty may not be uniform and prediction errors on tasks more critical to the optimization model tend to have a greater impact on generated solutions. Given the separation between the prediction and optimization phases, the prediction model cannot directly account for this interplay, potentially overlooking opportunities to fully harness predictive insights for enhancing optimization decisions.

Decision-aware prediction optimization (DAO) aims to overcome the separation of prediction and optimization by adapting the predictive modeling process to align closely with the specific optimization objective. DAO focuses on tailoring the predictive task to minimize \textit{decision loss}, defined as the negative impact of prediction inaccuracies on subsequent optimization efforts. Research in this field, such as \cite{elmachtoub_smart_2022, kallus_stochastic_2023} and \cite{chung_decision-aware_2022}, has primarily investigated how the loss functions of machine learning models can be adjusted to minimize the downstream decision loss. However, loss function formulation is just one element of machine learning. Decision loss and the effectiveness of optimization decisions can also be significantly impacted by other critical aspects of the machine learning to optimization pipeline, such as feature engineering, the choice of modeling algorithms, and the strategies employed for data aggregation and preparation. These considerations can even be significantly more important than the loss function in the learning process, particularly in scenarios characterized by complex prediction tasks and limited data.

To illustrate this, consider the task of learning individual worker behavior within a workforce. To learn each individual's unique behavior, one could train an individualized predictive model for each worker built using only that worker's data. These individualized models reflect each worker's unique decision-making history; however, as is common for human-centric processes, there may not be enough data to build a sufficiently complex model that represents the behavior of each worker. Alternatively, one could pool data across individuals to create an aggregate predictive model. This aggregated model benefits from a richer dataset, allowing for the application of more sophisticated modeling techniques and potentially enhancing overall performance. However, this approach overlooks the distinct nuances of individual workers' behaviors. This dilemma highlights a balance between achieving broad model performance and maintaining a degree of personalization. It is unclear if the choice that minimizes prediction loss and model accuracy aligns with minimizing decision loss within the downstream optimization problem.

In this paper, we propose and examine a framework that extends DAO to account for predictive model selection.  This enables an exploration of diverse predictive modeling methodologies, moving beyond the traditional focus on loss function formulation. Our framework considers a workforce comprised of individual workers, each with unique decision-making characteristics, and optimally selects a predictive model for each worker to be used in a downstream optimization problem. By integrating model selection into the DAO, we bridge the gap between individual worker tendencies and overarching organizational goals, enhancing the efficiency and effectiveness of decision-making processes. This project, conducted in collaboration with a U.S.-based auto-insurance provider, uses real-world data to demonstrate our framework's practical potential. The insurer employs underwriters to evaluate high-risk policy requests, with tasks currently assigned randomly across a large workforce. The insurer's leadership suspects that underwriters' decisions vary based on personal experience and risk tolerance. By modeling these individual tendencies using machine learning, we assess their impact within a downstream optimization framework, demonstrating that targeted predictive model selection can refine workforce allocation and create additional value for the organization.

The key contributions of this work, driven by the need for decision-aware model selection in real-world settings, are summarized as follows:
\begin{enumerate}
    \item \textbf{Decision-Aware Workforce Allocation Framework:} We propose a novel optimization framework that directly integrates predictive model selection into the workforce allocation process. Unlike traditional methods that treat prediction and optimization as separate steps, our DAO framework enhances the quality of worker-task assignments by dynamically selecting the most appropriate predictive model for each worker within the optimization process. This integration allows for a more accurate alignment of worker behavior with organizational objectives, resulting in improved decision-making and task allocation.
    \item \textbf{Predictive Modeling Approaches for Workforce Allocation:} We conduct an evaluation of three distinct predictive models to represent worker behavior in the context of workforce allocation. This analysis highlights the trade-offs between capturing the unique decision-making processes of individual workers and addressing the limitations imposed by data availability. By applying these models to real-world data, we demonstrate their effectiveness in optimizing worker-task assignments and explore their performance within both traditional UAO and our proposed DAO frameworks.
    \item \textbf{Application to Insurance Workforce Management:}  In collaboration with a leading U.S. auto-insurance provider, we apply our framework to real operational data and conduct an analysis of the results. This analysis includes evaluating the effectiveness of our approach in optimizing workforce allocation and discussing the managerial insights gained from these findings. Our study validates the practical applicability of our methods and provides actionable strategies for managers to enhance workforce management in real-world settings.
\end{enumerate}

The paper begins with a comprehensive literature review in Section \ref{sec:lit_review}. Section \ref{gen_form} introduces the workforce allocation problem and discusses methods for incorporating worker behavior into worker-task assignments, starting with traditional DPO and UAO approaches and then our proposed DAO framework. Section \ref{case_study_intro} presents a case study conducted with a U.S. auto-insurance provider to demonstrate the application of the DAO framework, covering the business context, predictive modeling approaches, and the implementation of UAO and DAO frameworks. Section \ref{sec_results} presents the findings of the case study and discusses key managerial insights. Finally, concluding remarks are provided in Section \ref{sec_conclusion}.


\section{Literature Review} \label{sec:lit_review}
Our research is related to four streams of literature: workforce allocation, predicting worker behavior, decision-aware optimization frameworks, and operations research applications in the insurance industry. 

\subsection{Workforce Allocation}
Assigning workers to tasks is a well-studied problem within operations research and has applications in many fields including production \citep{eiselt_employee_2008,vila_branch-and-bound_2014}, healthcare \citep{lanzarone_robust_2014}, and product management \citep{hendriks_human_1999}. We refer the reader to the review papers by \cite{van_den_bergh_personnel_2013} and \cite{bouajaja_survey_2017} for a comprehensive overview. 

Many authors study assignment problems while considering worker attributes such as preferences \citep{dowsland_nurse_1998,akbari_scheduling_2013}, productivity \citep{brusco_staffing_1998,brusco_exact_2008}, skills \citep{de_bruecker_workforce_2015}, and performance \citep{derman_sequential_1972,nikolaev_stochastic_2010}. Most similar to our work are those that individualize these attributes for the specific worker-task assignment. For example, \cite{dowsland_nurse_1998} consider the preferences of individual nurses for specific shifts, \cite{hanne_multiobjective_2009} explores company preference for worker-task assignments, and \cite{brusco_exact_2008} leverages departmental productivity of workers. However, these papers often assume worker attributes are known with certainty. Our work contributes to this stream of research by using data-driven methods to derive and understand worker attributes in subjective decision-making tasks, addressing the gap in the literature where uncertainty in worker attributes has been previously identified but not extensively explored \citep{de_bruecker_workforce_2015}. Furthermore, we explore the impact of uncertainty in these predictions within a downstream optimization problem.

\subsection{Predicting Worker Behavior}
The application of machine-learning techniques to understand and predict worker behavior is a growing field in operations research and human resource management. This interdisciplinary approach has gained momentum as organizations increasingly recognize the value of data-driven insights into workforce dynamics. Existing research includes predicting workforce attrition \citep{alduayj_predicting_2018, alsheref_automated_2022, raza_predicting_2022}, productivity \citep{florez_using_2017, alaloul_productivity_2022}, and performance \citep{chalfin_productivity_2016,sajjadiani_using_2019}. These papers typically focus on the methods used to identify worker attributes and the insights learned. Our paper moves beyond prediction and focuses on how the predictions can be used to improve decision processes within an organization. Most similar to our work is the paper by \cite{chalfin_productivity_2016} that discusses how individuals within a heterogeneous workforce should be hired or selected for promotions. However, unlike this work, we focus on how the method used to identify worker attributes impacts the value of the generated solution and resulting business decisions.

\subsection{Decision-Aware Learning}
Recent advancements in research are increasingly highlighting the importance of integrating predictive modeling with decision-making processes. This has led to the development of decision-aware approaches, which are specifically designed to align machine learning algorithms with the unique requirements of optimization problems. Notable contributions in this area include the work of \cite{elmachtoub_smart_2022} and \cite{kallus_stochastic_2023}, who have incorporated decision-loss functions into their models and adjusted decision forest algorithms based on their optimization effects, respectively. However, these approaches often encounter limitations in terms of their broader applicability, since they are heavily dependent on the specific structures of the optimization problems they address. To extend the reach of these methodologies, \cite{chung_decision-aware_2022} proposed a more flexible, model-agnostic approach, utilizing sample weighting derived from mathematical principles to increase their compatibility with a variety of machine learning models. Our study aligns with the approach of \citeauthor{chung_decision-aware_2022} in maintaining the conventional structure of machine-learning training algorithms. However, we differentiate our research from these approaches by focusing on the selection of predictive models instead of modifying the training algorithms themselves. Similarly, while \cite{wilder_melding_2019} have explored similar decision-aware approaches in the realm of combinatorial optimization, their emphasis on machine-learning training algorithms contrasts with our approach of predictive model selection. Our research contributes to the field of decision-aware methodologies by advocating for a more comprehensive strategy that not only considers the training of predictive models but also emphasizes the importance of strategic model selection to optimize decision-making processes.

\subsection{Operations Research in the Insurance Industry}
Lastly, this work contributes to the stream of research that applies operations research methods in the insurance industry. Contributions in this field have focused on using techniques such as game theory, dynamic programming, and stochastic models to address challenges related to underwriting, risk management, and claims processing \citep{brockett_operations_1997, dedu_multiobjective_2015, consigli_dynamic_2011, wuthrich_machine_2018}. In particular, recent advancements have increasingly incorporated machine learning techniques to improve predictive accuracy and operational efficiency in insurance, focusing on pricing optimization, fraud detection, and risk classification \citep{blier-wong_machine_2020, santos_machine_2024}. Our work aligns most closely with the study by \cite{bertsimas_algorithmic_2022}, who explore how machine-learning classification thresholds impact risk management in algorithmic insurance. However, rather than treating classification thresholds as parameters we can easily control, we focus on accounting for the human aspect of decision-making, particularly the variability in human behavior and judgment, as it influences the predictive modeling and optimization of worker-task assignments in insurance underwriting.

\section{Formulation} \label{gen_form}
In this section, we introduce the workforce allocation problem and outline methods for leveraging worker behavior when making worker-task assignments. In Section~\ref{workforce_allocation}, we begin by presenting a general workforce allocation problem and discuss traditional approaches for incorporating predictive insights into worker-task assignments. Although these methods are straightforward, they treat prediction and optimization as separate steps, neglecting how the predictive model used to represent worker behavior influences the resulting assignment policy. To address this disconnect and improve solution quality, in Section~\ref{sec_DAO}, we propose a novel decision-aware workforce allocation framework. This method enhances worker assignments by directly integrating predictive model selection for worker behavior into the optimization framework.

\subsection{Workforce Allocation Problem} \label{workforce_allocation}
We first introduce the general workforce allocation problem. In this problem, a workforce of various workers is allocated to various tasks. Let $I$ denote the set of workers and $J$ the set of tasks. We define the decision matrix $X \in \{0,1\}^{|I| \times |J|}$, where $x_{ij} = 1$ if worker $i$ is assigned to task $j$, and $x_{ij} = 0$ otherwise. For now, we keep the remaining feasible region of assignments, $\mathcal{X}$, general but note that worker assignments are typically constrained by factors such as worker capacity and skill prerequisites.

Once a worker is assigned tasks, each task requires a binary decision from the assigned worker. For example, such tasks might include deciding whether to accept or reject an auto insurance request, prescribing medication to a patient, or making a hiring decision for a job candidate. We introduce the Bernoulli random variable $B_{ij} \in \{0,1\}$ to represent the behavior of worker $i$ when assigned to task $j$. In particular, $B_{ij} = 1$ indicates that worker $i$ takes an action (e.g., approving a request) on task $j$, while $B_{ij} = 0$ indicates that an action is not taken. The value of performing an action for task $j$ is denoted by $v_j \in \mathbb{R}$, where $v_j$ can be positive or negative, depending on the associated risks or benefits of the task. Our objective is to maximize the total expected value of the actions taken, subject to the constraints defined by the feasible region $\mathcal{X}$. Formally, the workforce allocation problem is defined as:
\begin{subequations}
\begin{align}
\text{maximize } & \sum_{i \in I} \sum_{j \in J} \mathbb{E}\left[B_{ij}\right] v_j x_{ij} \label{obj_fun_stochastic}\\
\text{subject to } &\  
x_{ij} \in \{0,1 \}, \quad \forall i \in I, j \in J,
\\ 
&X \in \mathcal{X}.
\end{align}
\end{subequations}

The primary challenge in this optimization model is the treatment of the expected value in the objective. One naive approach is to assume we have knowledge of worker behavior and replace $E[B_{ij}]$ with a deterministic parameter. In reality, due to the subjective nature of tasks, $B_{ij}$ is often uncertain until worker behavior is observed. To inform worker behavior, a reasonable approach is to utilize a predictive modeling technique such as a machine-learning model. These predictive models can utilize historical data to predict worker actions, which can be used in optimizing worker allocation.

\subsubsection{Decoupled Prediction Optimization (DPO).} \label{sec_DPO}
The simplest method to integrate predictive modeling into an optimization framework is the DPO approach. In this method, a trained predictive model $m$ is used to directly predict the action of worker $i$ on task $j$. Specifically, let $B^{m}_{ij} \in \{0,1\}$ denote the predicted action made by model $m$ for worker $i$ on task $j$. The DPO formulation is as follows:
\begin{subequations}
\begin{align}
\text{maximize } &\ \sum_{i \in I} \sum_{j \in J} B^{m}_{ij} v_j x_{ij} \label{obj_fun_dpo}\\
\text{subject to } &\  
x_{ij} \in \{0,1 \}, \quad \forall i \in I, j \in J,
\\ 
&X \in \mathcal{X}.
\end{align}
\end{subequations}
Although the decision maker typically selects a predictive model that maximizes prediction accuracy, a shortcoming of this approach is the assumption that these predictions accurately represent $B_{ij}$. In practice, it can be difficult to assess if the optimization model is exploiting the true behavior of each worker or potential inaccuracies of the predictive models. If the latter, then the optimization model may produce solutions that perform poorly in practice.

\subsubsection{Uncertainty-Aware Optimization (UAO).} \label{sec_UAO}
While the DPO approach offers a straightforward framework for integrating predictive modeling with optimization tasks, it assumes that the predictions are equally accurate. To mitigate this assumption, UAO is based on the idea that $\mathbb{P}(B_{ij}^{m} = 1) \approx \mathbb{E}[B_{ij}]$. This representation intuitively captures the models confidence that $B_{ij}=1$, which may vary across $i-j$ pairs. The UAO formulation can be written as:
\begin{subequations}
\begin{align}
\text{maximize } &\ \sum_{i \in I} \sum_{j \in J} \mathbb{P}(B^{m}_{ij} = 1)v_j x_{ij} \label{obj_fun_uao}\\
\text{subject to } &
x_{ij} \in \{0,1 \}, \quad \forall i \in I, j \in J,
\\ 
&X \in \mathcal{X}.
\end{align}
\end{subequations}
Although the UAO approach captures predictive model confidence, $\mathbb{P}(B^{m}_{ij} = 1)$ can vary depending the choice of predictive model (e.g., decision tree versus logistic regression), data preparation techniques, and training process. In particular, for UAO, predictive model selection occurs independently from the optimization process.

\subsection{Decision-Aware Optimization (DAO)} \label{sec_DAO}
For many workforce settings, predicting human behavior is challenging due to the complexity of decision-making processes and the limited data typically generated from manual human activities. Traditional machine learning approaches gauge the effectiveness of predictive modeling using metrics such as accuracy or area-under-the-curve (AUC). However, these metrics may not fully capture the practical implications of prediction errors in operational contexts.

Our proposed DAO framework addresses this challenge by embedding predictive model selection directly within the optimization process. Unlike DPO and UAO approaches where the outputs of predictive models are treated as input to the optimization model, DAO acknowledges that predictions influence the decisions made within the optimization framework. Our framework aims to adapt predictive model selection to better support decision making, thereby enhancing the overall effectiveness and reliability of the optimization model solution.

To implement the DAO approach, we consider a set of candidate predictive models, denoted by $ M(i) $, for predicting each worker's behavior. Each model $ m \in M(i) $ estimates the behavior of worker $ i $ on task $ j $, represented by $ B_{ij}^m  \in \{ 0,1\}$. We introduce binary decision variables $ y_{im} $ to indicate the predictive model selection, where $ y_{im} = 1 $ indicates that model $ m $ is selected to represent worker $ i $ and $ y_{im} = 0 $ if predictive model $ m $ is not selected for worker $i$. Additionally, we introduce binary variables \( z_{ijm} \) to link the work allocation with predictive model selection. The variable \( z_{ijm} = 1 \) indicates that worker  \( i \) is allocated to task \( j \) using the predictions from model \( m \) and \( z_{ijm} = 0 \) otherwise. The DAO formulation can be written as:

\begin{subequations}\label{eq:DAO}
\begin{align}
\text{maximize } &\ \sum_{i \in I} \sum_{j \in J} \sum_{m \in M(i)} \mathbb{P}(B^{m}_{ij} = 1) v_j z_{ijm} \label{obj_fun_dao} \\
\text{subject to }
& \sum_{m \in M(i)} z_{ijm} = x_{ij}, \quad \forall i \in I, j \in J, \label{cons1} \\
& z_{ijm} \le y_{im}, \quad \forall i \in I, j \in J, m \in M(i) \label{cons2} \\
& \sum_{m \in M(i)} y_{im} = 1, \quad \forall i \in I, j \in J, \label{y_cons} \\
& z_{ijm} \in \{ 0, 1\}, \quad \forall i \in I, j \in J, m \in M(i), \label{binary_cons3} \\
& y_{im} \in \{ 0, 1\}, \quad \forall i \in I, m \in M(i), \label{binary_cons2} \\
& x_{ij} \in \{ 0, 1\}, \quad \forall i \in I, j \in J, \label{binary_cons} \\
&X \in \mathcal{X}. \label{X_cons} 
\end{align}
\end{subequations}
The objective function \eqref{obj_fun_dao} aims to maximize the expected value of workforce allocations by simultaneously optimizing predictive model selection and worker-task assignments. Constraint sets \eqref{cons1} and \eqref{cons2} ensure each worker-task assignment is represented by one of the selected predictive models. Constraint set \eqref{y_cons} ensures that each worker is assigned exactly one predictive model from their set \( M(i) \). Constraint sets \eqref{binary_cons3}-\eqref{binary_cons} enforce the binary requirements on the decision variables. Constraint \eqref{X_cons} enforces the worker allocation requirements represented by the set $\mathcal{X}$. 

We note that the DAO formulation generalizes both the UAO and DPO approaches. The DAO formulation simplifies to a UAO model when \( |M(i)| = 1 \forall i \in I\). In this case, the predictive model used to represent each worker is predetermined (i.e., $y_{im} = 1$) and $x_{ij} = z_{ijm}$. The DAO formulation simplifies to the DPO approach when \( |M(i)| = 1 \) and the $B_{ij}^m$ are assumed to be deterministic (i.e.,$\mathbb{P}(B^{m}_{ij} = 1) = B^{m}_{ij}$). Lastly, we note that the constraint set \eqref{y_cons} serves as a form of regularization because it restricts the selection of one predictive model per worker, which helps to prevent overfitting. While relaxing this constraint could introduce more flexibility in model selection, we retain it to ensure simplicity, interpretability, and to reduce the risk of overfitting by avoiding overly complex model choices for individual workers.

%
\subsubsection{Adjusting Model Predictions.} \label{sec_rep_preds}
A challenge when assigning tasks to workers using multiple predictive models is that the accuracy of each model \( m \in M(i) \) can vary significantly depending on the task. For example, if a predictive model \( m \in M(i) \) is trained using data from multiple workers, it may perform well for tasks where there is a high level of consensus among workers. However, for tasks where decision-making varies significantly across workers, the same model might perform poorly, as it may fail to capture the subtleties of individual worker behavior in more subjective scenarios.

Furthermore, in subjective decision-making contexts, tasks are often unique, making it difficult for predictive models trained on historical data to generalize effectively to new tasks. For instance, a random forest model might estimate the probability \( \mathbb{P}(B_{ij}^{m} = 1) \) as the proportion of decision trees that predict \( B_{ij} = 1 \), based on past worker behavior.  However, given the diverse nature of tasks that require nuanced human judgment, new tasks may introduce variations not previously captured in the training data. In such settings, relying solely on historical training data can lead to predictions biased by the observed training data and may not accurately reflect the true dynamics of current tasks. 

To mitigate this issue, we propose adjusting the expected behavior prediction \( \mathbb{E}[B_{ij}] \) by incorporating the predictive model's performance on a test dataset. This adjustment helps account for the model's ability to generalize beyond the training data and ensures a more accurate representation of worker behavior for future tasks.

We begin by expressing the expected behavior \( \mathbb{E}[B_{ij}] \) as follows:
\begin{subequations}
\begin{align}
\mathbb{E}[B_{ij}] &= \mathbb{P}(B_{ij}^{m} = 1),\\
&= \mathbb{P}(B_{ij}=1 \mid B_{ij}^{m} = 1) \cdot \mathbb{P}(B_{ij}^{m} = 1) \\
& \quad + \mathbb{P}(B_{ij}=1 \mid B_{ij}^{m} = 0) \cdot \mathbb{P}(B_{ij}^{m} = 0). \notag
\end{align}
\end{subequations}

In this equation \( \mathbb{P}(B_{ij}^{m} = 1) \) and \( \mathbb{P}(B_{ij}^{m} = 0) \) represent the likelihood that the predictive model \( m \) forecasts a specific action by worker \( i \) on task \( j \). The terms \( \mathbb{P}(B_{ij} = 1 \mid B_{ij}^{m} = 1) \) and \( \mathbb{P}(B_{ij} = 1 \mid B_{ij}^{m} = 0) \) capture the conditional probabilities that worker \( i \) actually takes the action, given the model’s prediction. Since these conditional probabilities \( \mathbb{P}(B_{ij} = 1 \mid B_{ij}^{m} = 1) \) and \( \mathbb{P}(B_{ij} = 1 \mid B_{ij}^{m} = 0) \) cannot be directly observed until worker behavior is observed, we estimate these terms using the model’s performance on a test dataset. This allows us to better gauge how well the model generalizes to new tasks. Thus, we adjust the expected behavior prediction as follows:
\begin{align}
\mathbb{E}[B_{ij}] 
&\approx \hat{\mathbb{P}}(B_{ij}=1 \mid B_{ij}^{m} = 1) \cdot \mathbb{P}(B_{ij}^{m} = 1) \label{E_approx}\\
& \quad + \hat{\mathbb{P}}(B_{ij}=1 \mid B_{ij}^{m} = 0) \cdot \mathbb{P}(B_{ij}^{m} = 0). \notag
\end{align}
Note that \( \mathbb{P}(B_{ij}^{m} = 1) \) and \( \mathbb{P}(B_{ij}^{m} = 0) \) are known from the predictive model, while the terms \( \hat{\mathbb{P}}(B_{ij} = 1 \mid B_{ij}^{m} = 1) \) and \( \hat{\mathbb{P}}(B_{ij} = 1 \mid B_{ij}^{m} = 0) \) represent the estimated conditional probabilities of the worker’s actual behavior, derived from the test dataset's results. Several techniques can be employed to estimate these conditional probabilities, including clustering methods and simple predictive models. Furthermore, when the prediction error is uniformly distributed, \( \hat{\mathbb{P}}(B_{ij} = 1 \mid B_{ij}^{m} = 1) \) and \( \hat{\mathbb{P}}(B_{ij} = 1 \mid B_{ij}^{m} = 0) \) correspond to the positive predictive value and the false omission rate of model \( m \in M(i) \) when applied to worker \( i \) \citep{parikh_understanding_2008}.

By adjusting model predictions based on test dataset performance, we integrate out-of-sample predictive performance into our decision-making process. This method enhances both the accuracy and relevance of task assignments. Additionally, employing estimated conditional probabilities offers a richer analysis than straightforward accuracy metrics alone. This approach aligns model predictions with individual decision-making tendencies, including factors like risk aversion and tolerance, ensuring that assignments reflect true worker behavior more closely.


\section{Case Study: Insurance Underwriter Allocation} \label{case_study_intro}

To demonstrate the application of the proposed DAO framework, we partnered with a leading U.S.-based auto-insurance provider that manages a workforce of insurance underwriters who are responsible for reviewing and making decisions on unusual or high-risk insurance policy requests. The insurance industry offers a high-stakes financial environment, making it an ideal setting for exploring workforce allocation challenges.

In this section, we present the details of this case study. First, in Section~\ref{sec_bc}, we provide the motivation and details of the business context. Next, in Section~\ref{sec_learn_behavior}, we present the real-world dataset and propose three predictive modeling approaches for learning underwriter behavior. In Section~\ref{sec_optimization}, we introduce the UAO and DAO models used to allocate workers to tasks. The results of the case study, along with key managerial insights, are presented in Section~\ref{sec_results}.

\subsection{Business Context} \label{sec_bc}
In this partnership, we examine the allocation of insurance underwriters. Underwriters are employees that evaluate the details and terms of high-risk or unusual insurance policy requests and determine whether the company should approve the requested coverage. This decision-making process involves a binary choice, either to accept or reject the policy. High-risk policies have the potential to be particularly profitable for the insurance company, as customers typically pay a higher premium (i.e., the payment made by the customer for insurance coverage). However, these policies also come with a greater risk of financial loss due to the inherent uncertainty and unpredictability associated with insuring high-risk individuals or situations. 

Traditionally, the task of reviewing a particular high-risk policy is assigned to an underwriter randomly. Although underwriters undergo standardized training, we hypothesize that individual differences, such as varying levels of experience and risk tolerance, may lead to variations in their performance. Our objective is to identify these tendencies and demonstrate how incorporating such insights into the workforce allocation process can support an organizational goal, without creating additional work for the employees. For example, if the insurance company believes a policy should be accepted, it may wish to avoid assigning it to an underwriter who has historically shown aversion to similar policies.


To inform worker allocation, we need to predict how underwriters will behave on future policies, which is challenging for two key reasons. First, each underwriter deals with a wide variety of unique policy requests, each with distinct characteristics. This variability makes it difficult to generalize patterns across policy requests as the decision criteria for one type of policy may not apply to another. Second, because underwriting is a manual and judgment-based process, each underwriter generates only a limited amount of data that captures their decision-making patterns. This scarcity of data further complicates predictive modeling as the available data may not fully reflect the nuances of each underwriter's approach to different policy types. These limitations underscore the need for the strategic application of predictive modeling insights. This is where the DAO method proves valuable, enabling more tailored and effective assignments that align with underwriter strengths.

\subsection{Learning Underwriter Behavior} \label{sec_learn_behavior}
In learning underwriter behavior, our goal is to predict whether a particular underwriter will accept or reject a policy given the details of the policy request. The dataset for this case study was provided by the auto-insurance company and includes deidentified requests for auto insurance coverage along with the decisions made by the assigned underwriters. The policy details include relevant information such as reasons for underwriter review, customer demographics (age and years driving), vehicle specifications, intended vehicle use, and requested coverage. After removing correlated features (above $0.85$ or below $-0.85$) and before converting categorical variables to binary features (i.e., one-hot encoding), the dataset contains 56 features in total.

We use 15 months of data to develop our predictive models. Twelve months of data were used for training the predictive models, and data from the subsequent 3 months were used for testing. A final month of data is reserved for the optimization phase of the case study. Underwriters that reviewed fewer than 50 policies in the testing dataset were removed and any policy reviewed by more than one underwriter was excluded. After applying these filters, we were left with data for 53 underwriters; 17,853 rows for training and 3,857 rows for testing. The data inclusion process is outlined in Appendix A.

Using this data set, we develop three approaches to predict if an underwriter approves a given policy based on the features of that policy. The approaches include constructing individualized models for each underwriter (Section \ref{sec_ind}), developing an aggregate model that leverages all available data (Section \ref{sec_agg}), and employing a profile-based approach that combines elements of both individual and aggregate modeling (Section \ref{sec_prof}). This range of strategies provides the flexibility to select a predictive model that best aligns with each underwriter's unique behavior and the data available.

\subsubsection{Individual Predictive Model.} \label{sec_ind} 
We first consider individualized predictive models. For this approach, we develop a unique predictive model for each worker, trained on data generated by that worker's previous decisions. While this method captures each worker's unique decision-making processes and behavioral tendencies, it can be limited by the availability of sufficient data, an issue commonly encountered with new hires or in rapidly evolving industries.

We use a logistic regression model with L2 (ridge) regularization for each worker. Logistic regression is selected because it performs well in settings with limited training data and L2 regularization provides robustness against outliers in data sets with many features, particularly in data-scarce settings. The weight of the regularization term is tuned using 10-fold cross-validation on the training set. We use area-under the receiver operating characteristic curve (AUC) as our primary metric.

Figure~\ref{fig:Ind_Agg} shows the performance of each individual predictive model quantified by  test set AUC. An AUC of 0.5 is equivalent to guessing and larger values indicate higher predictive power. The weighted (by number of observations) average AUC across all individual models was 0.62. Notably, only 11.3-percent of the models achieved an AUC above 0.70, suggesting that the majority of individual predictive models perform only slightly better than randomly guessing a worker's decisions \citep{hanley_meaning_1982,hosmer_jr_applied_2013}.

The low AUC values of the individual predictive models indicate that while these models can provide personalized insights, their overall predictive performance is modest. We hypothesize the low AUC values are primarily due to the limited training data available per worker. In the workforce allocation context, these results highlight an important trade-off. While individual models created using personalized data may capture personalized worker behavior, the effectiveness of these models can be limited by data scarcity. Therefore, in data-limited settings, relying solely on individual predictive models may not yield strong performance. This finding highlights the need to incorporate other predictive modeling strategies, such as aggregate or profile-based models, to improve predictive accuracy.

\subsubsection{Aggregate Predictive Model.} \label{sec_agg}
To address data scarcity issues that may impact individual predictive model performance, we propose an aggregate predictive model. In this approach, a single predictive model is trained using the available data from all workers with the same features as the individual predictive models. For tasks where there is a consensus among workers, this approach is reasonable as it provides valuable insights into general decision-making patterns. However, a potential limitation is its lack of sensitivity to individual differences, leading to inaccuracies in settings where consensus is low. 

Unlike in the individual predictive model, where we are constrained by the need for models that perform well with smaller datasets, the aggregated data allows us to train a random forest model, capable of capturing complex non-linear relationships due to the larger volume of data available. The minimum number of samples in each leaf of the random forest model is tuned using 10-fold cross-validation on the training set.

\begin{figure}[ht!]
    \centering
    \includegraphics[width = .75\textwidth]{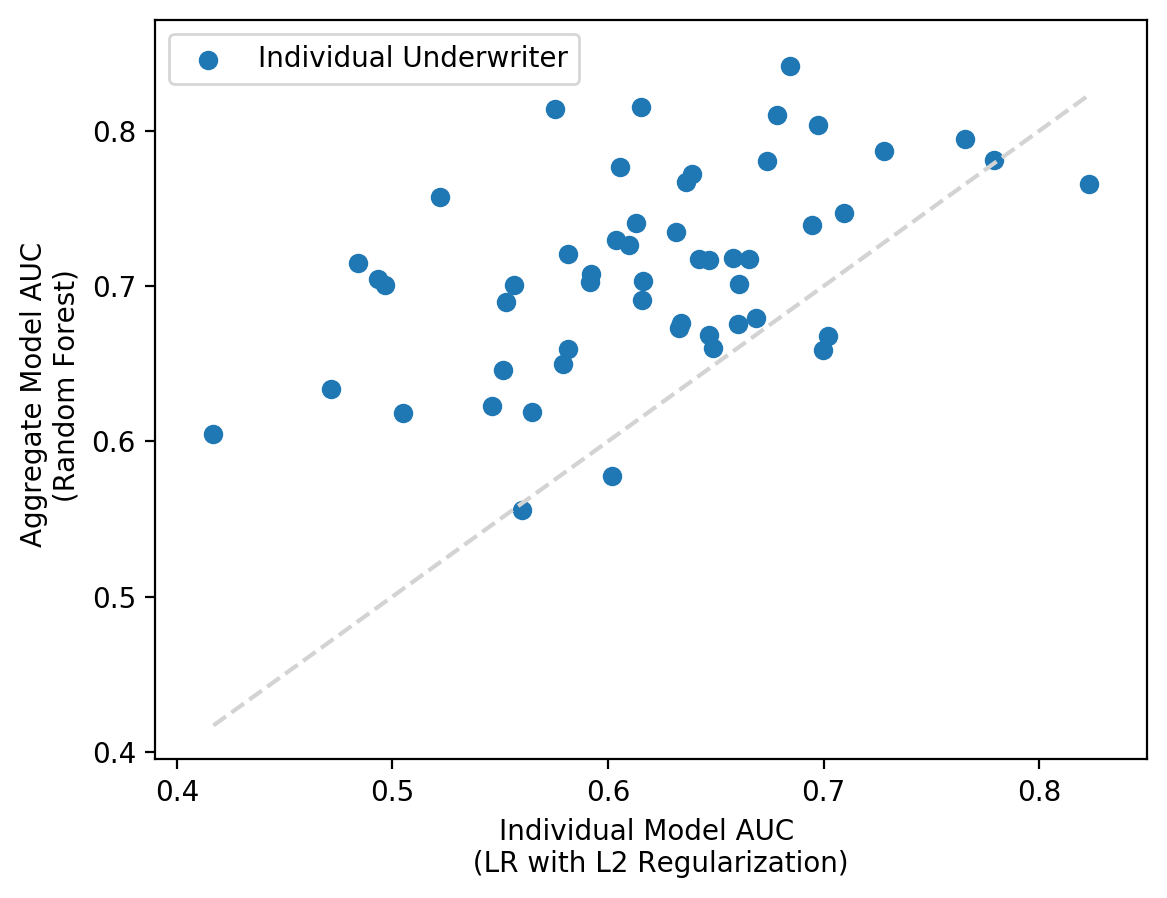}
    \caption{Comparison of the test set AUC between individual and aggregate predictive models for each underwriter, with the grey dashed line indicating where the individual and aggregate models achieve equal AUC values.} \label{fig:Ind_Agg}
\end{figure}

Figure \ref{fig:Ind_Agg} displays the AUC values obtained from applying the aggregate model to each individual underwriter in the test set. The majority of points lie above the 45-degree line, indicating that the aggregate predictive model achieves a better test set AUC value than the individual predictive model for most underwriters. In fact, the weighted (by number of observations) average AUC of the aggregate predictive model across all individual underwriters was 0.72; an increase of 16-percent over the individual models. Additionally, when applied to specific workers, the aggregate model achieved an AUC of 0.70 or higher for 62.2-percent of the workers, suggesting that it provides moderate to good predictive power for the majority of underwriters \citep{hanley_meaning_1982,hosmer_jr_applied_2013}.

\begin{figure}[ht!]
    \centering
    \includegraphics[width = .8\textwidth]{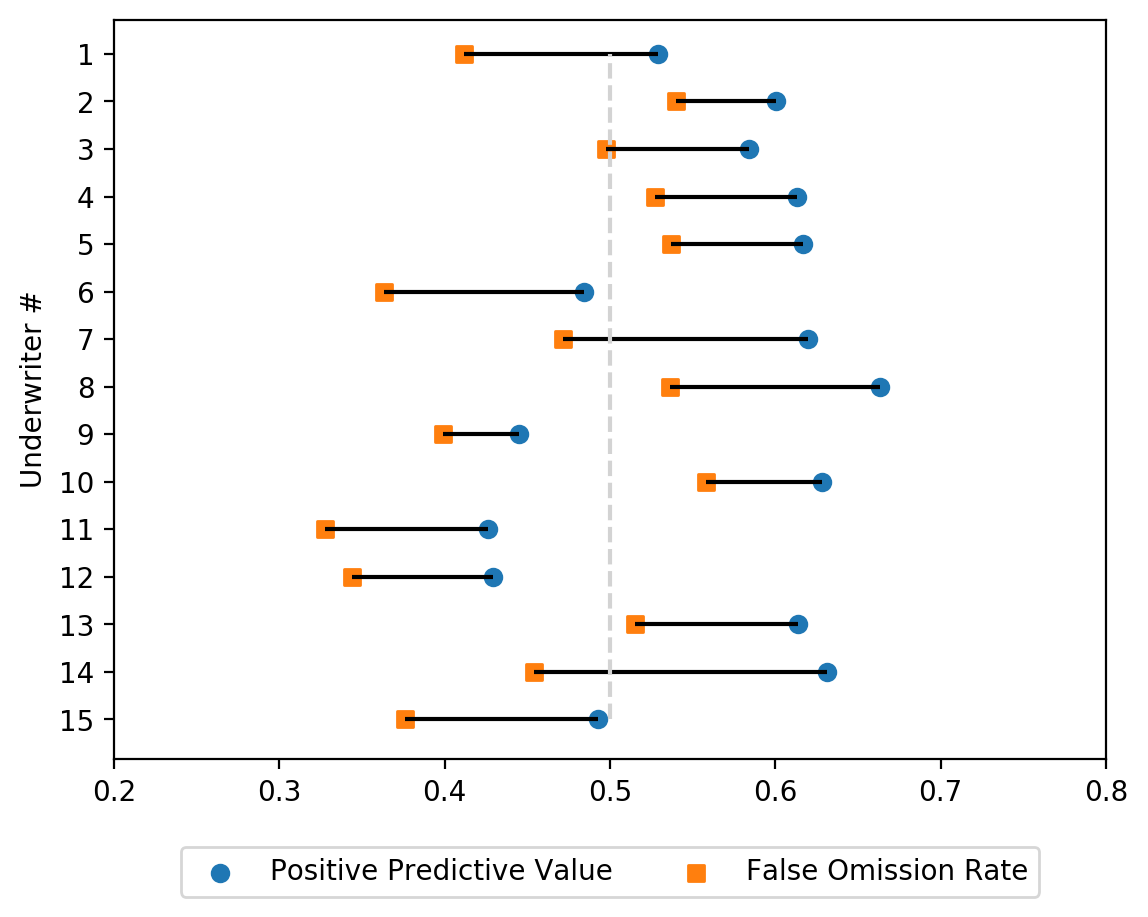}
    \caption{Out-of-sample Positive Predictive Value (PPV) and False Omission Rate (FOR) of the aggregate predictive model when applied to fifteen underwriters.} \label{fig:Gen_Performance}
\end{figure}

To further evaluate the aggregate model’s fit for each worker, we compare each worker's actual behavior to the aggregate model's predictions. Using the test dataset, we analyze the Positive Predictive Value (PPV) and False Omission Rate (FOR) of the aggregate model for each worker. The results are summarized for 15 of the 53 underwriters with the most training data in Figure \ref{fig:Gen_Performance}. A high PPV and low FOR (wide bar) indicate a strong fit. When both metrics are high, the underwriter is more risk-tolerant (bar skewed right), whereas when both are low (bar skewed left), the underwriter is more risk-averse. As shown in Figure \ref{fig:Gen_Performance}, workers vary in how well they are represented by the aggregate model and in how they behave with respect to it. This suggests that while the aggregate model may work well for some workers, it may not be suitable for others, further motivating the DAO approach.

While computing PPV and FOR assesses individual performance against the aggregate model, it does not capture individual tendencies that might be discerned through the more personalized approach of the individual predictive model. The aggregate model is designed to reflect broader trends and commonalities rather than the unique subtleties that might be evident in individual-specific models. Therefore, for worker-allocation, using an aggregate predictive model is more effective in settings where collective behavior outweighs individual differences. However, in environments where decision-making varies significantly across workers, a more personalized predictive modeling approach may be required.

\subsubsection{Profile Predictive Model.} \label{sec_prof}
To balance the trade-offs between individualized and aggregate predictive models, we also examine a profile-based predictive model. In this approach, workers are clustered according to distinct behavioral patterns, and a predictive model is built for each cluster. This method aims to mitigate the data limitations inherent in individual models while providing more nuanced predictions than the aggregate model.

\begin{figure}[ht!]
    \centering
    \includegraphics[width = .8\textwidth]{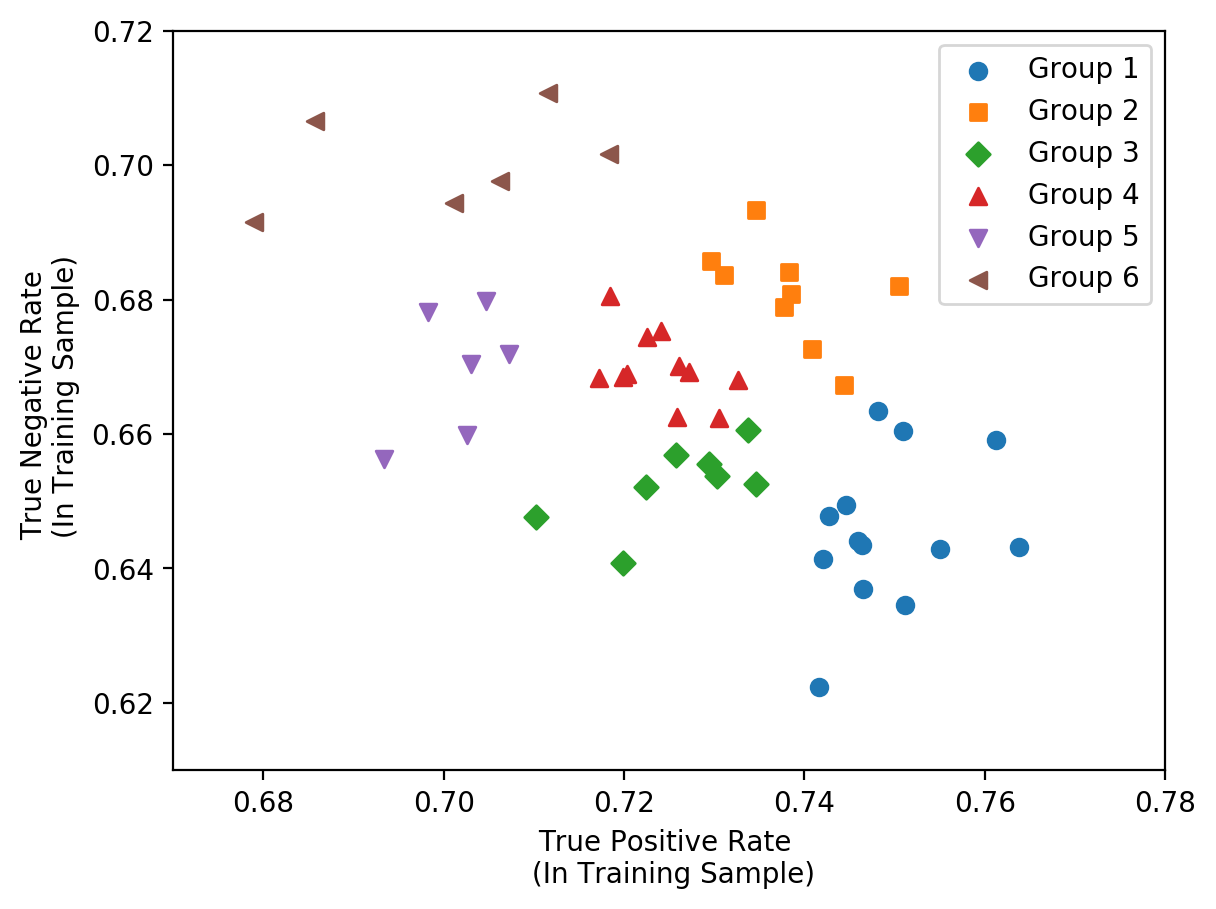}
    \caption{In-sample true positive and true negative rates are used to cluster workers for the creation of profile predictive models.} \label{fig:clusters}
\end{figure}

We clustered underwriters by assessing their actual behavior relative to the aggregate predictive model presented in Section \ref{sec_agg}. This categorization was achieved by calculating the true positive and true negative rates within the training sample. Subsequently, k-means clustering was applied to group workers based on these performance metrics. The results of this analysis are displayed in Figure \ref{fig:clusters}, where the true positive and true negative rates for each underwriter are plotted on the x-axis and y-axis, respectively. The clustering process resulted in six distinct clusters, as shown by the shapes and coloring in the figure. The decision to use six clusters was based on the distribution of workers and a thorough analysis of inter-cluster distances, which is further detailed in Appendix B. These clusters represent distinct behavioral groups within the workforce. The data generated by the underwriters in each group is then used to develop a profile predictive model. We employed logistic regression models with L2 regularization for each cluster, similar to the individual predictive models, to avoid issues that may arise from limited data. 

The weighted (by number of observations) average AUC of the profile predictive models across all individual underwriters was 0.67; higher than the individual models and lower than the aggregate.  When the profile models were applied to the specific underwriters whose data was used to generate them, 16.7-percent of the underwriters' models had an AUC below 0.70.


\subsubsection{Discussion.}
Depending on the specific context and the data at hand, each predictive modeling approach has its unique advantages. Individual models are particularly effective for in-depth analysis of a single worker's behavior, although their performance can be hindered by a scarcity of data. The aggregate model is better at capturing complex patterns consistent across workers, but may not adequately capture the subtle nuances of each individual's decision-making process. The profile predictive model blends these two approaches, integrating the detailed insights of individual behaviors with the comprehensive trends observed in aggregate data. 

\begin{figure}[ht!]
    \centering
    \includegraphics[width = .8\textwidth]{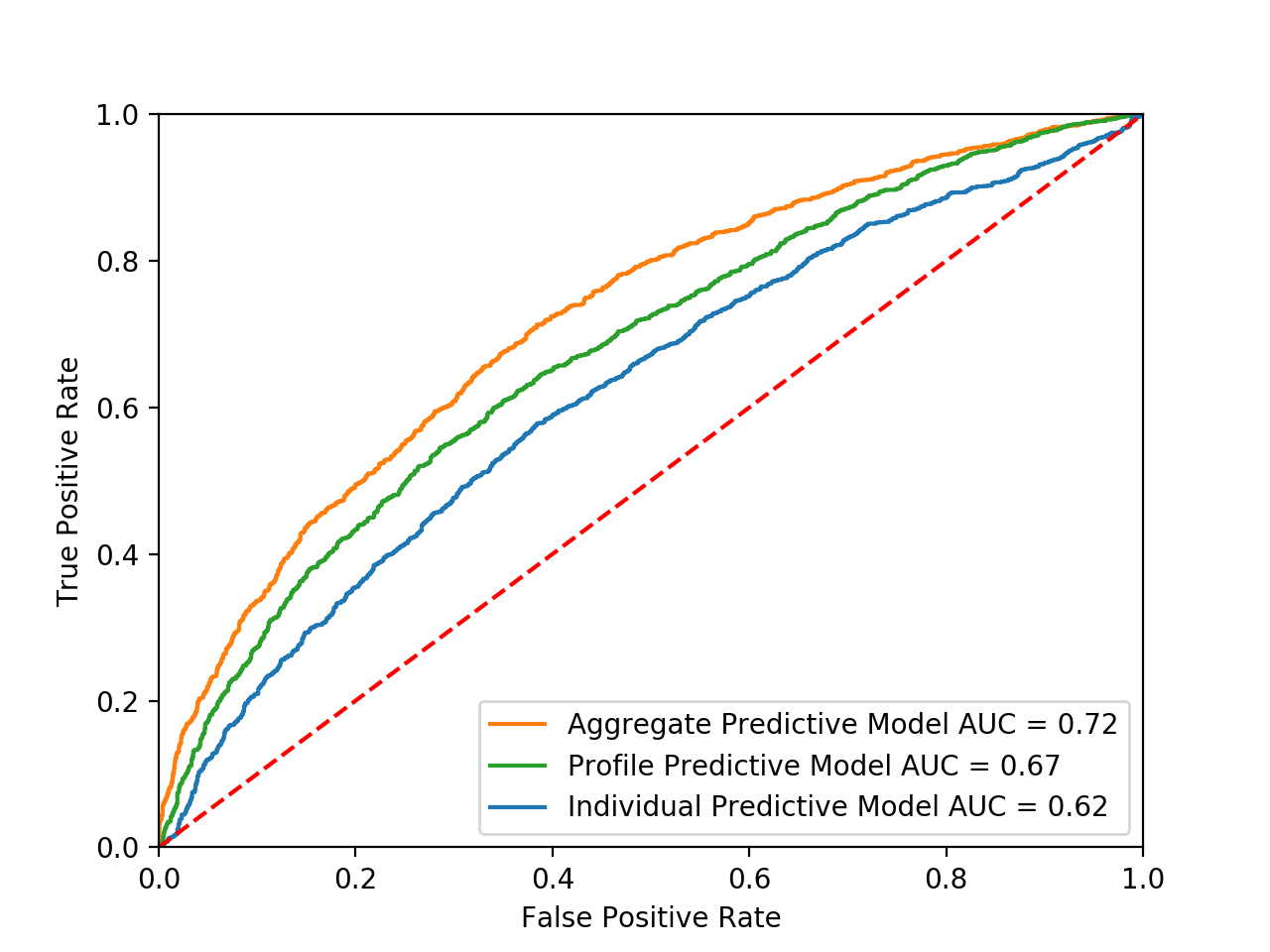}
    \caption{ROC curves and AUC values illustrating the predictive power of the proposed predictive models.} \label{fig:ROC}
\end{figure}

The effectiveness and comparative performance of these models can be visualized in Figure \ref{fig:ROC}, which presents a receiver operating characteristic (ROC) curve and AUC values, illustrating the predictive power of each model when applied to a test data set. As shown in Figure \ref{fig:ROC}, the aggregate predictive model has the highest AUC value of 0.72, indicating higher predictive power. The individual model shows the least predictive power with an AUC of 0.62, while the profile model is in-between with an AUC score of 0.67.

This analysis can be expanded to explore other predictive models, such as different worker groupings or more complex modeling techniques, to potentially enhance prediction accuracy further. The selected predictive models are meant to highlight a range of options that can be used in this type of setting. The choice of predictive model candidates ultimately depends on the specific goals of the analysis, the data available, and the behaviors studied.

In traditional DPO and UAO, we would typically select the predictive model with the highest AUC or apply the same approach uniformly across all workers. However, it is important to recognize that the predictive model with the highest predictive power, as indicated by AUC values, does not always guarantee the best decision performance. In the workforce allocation optimization model, the utility of each predictive model depends not only on its standalone accuracy but also on how each worker is utilized within the workforce allocation framework. In the next section, we explore this idea by applying these methods to the workforce allocation problem, using both the UAO and DAO approaches.

\subsection{Optimization Model} \label{sec_optimization}
We apply these predictive models to the proposed workforce allocation frameworks to demonstrate importance of predictive model selection. Initially, we assume that each worker's predictive model is predetermined and apply the UAO framework as outlined in Section \ref{sec_UAO}. This involves implementing the UAO framework using each of the three predictive modeling methods (individual, aggregate, and profile predictive models). Then, we use the DAO framework and optimize predictive model selection within the assignment model itself, as detailed in Section \ref{sec_DAO}. We evaluate the effectiveness of the optimized assignment strategies using the one month of data that was set aside from the training and testing datasets.

For the one-month planning horizon, the set $I$ represents the underwriters in the workforce, with $|I| = 53$. The set $J$ includes the policies that need to be reviewed, with $|J| = 1,187$. The set $M(i)$ contains the predictive models that represent the behavior of underwriter $i \in I$. In the UAO setting, $|M(i)| = 1$, meaning each underwriter has one predetermined predictive model. In the DAO setting, $|M(i)| = 3$, representing the choice among the three predictive models.

We compute the probability terms $\mathbb{E}[B_{ij}]$ using the approximation presented in equation \eqref{E_approx}. We estimate $\hat{\mathbb{P}}(B_{ij} = 1 \mid B_{ij}^{m} = 1)$ and $\hat{\mathbb{P}}(B_{ij} = 1 \mid B_{ij}^{m} = 0)$ as the PPV and FOR of model $m \in M(i)$, when applied to underwriter $i$, respectively. We assume that the parameters $v_j$ are known, representing the profit associated with accepting the policy under review in task $j$. We inform these parameters using a separate machine-learning model, noting that insurance company may have internal models capable of estimating the expected value of accepting various policies.

Finally, we define the set $\mathcal{X}$ to ensure that the assignments accurately reflect the time-based and capacity constraints of the real-world business context. Let $T$ represent the set of days within the one-month planning horizon, and $J_t$ denote the policies reviewed during day $t \in T$. To maintain alignment with historical demand, we require that each task be reviewed on the same day as it was in the historical data. Additionally, each underwriter is assigned a capacity, $c_{it}$, representing the number of tasks that worker $i \in I$ reviewed during day $t$ in the historical data. We assume that each underwriter cannot exceed this daily capacity. With these constraints defined, we formulate the complete DAO formulation as follows:

\begin{subequations}
\begin{align}
\text{maximize } &\ \sum_{i \in I} \sum_{j \in J} \sum_{m \in M(i)} \mathbb{P}(B^{m}_{ij} = 1) v_j z_{ijm} \label{eq:full_obj_fun_dao} \\
\text{subject to }
& \sum_{m \in M(i)} z_{ijm} = x_{ij}, \quad \forall i \in I, j \in J, \label{eq:full_cons1} \\
& z_{ijm} \le y_{ij}, \quad \forall i \in I, j \in J, m \in M(i) \label{eq:full_cons2} \\
& \sum_{m \in M(i)} y_{im} = 1, \quad \forall i \in I, j \in J, \label{eq:full_y_cons} \\
&\sum_{i \in I} x_{ij} = 1, \quad \forall j \in J, \label{eq:full_x_1}\\
    & \sum_{j \in J_t} x_{ij} \leq c_{it}, \quad \forall i \in I, t \in T, \label{eq:full_x_2} \\
& z_{ijm} \in \{ 0, 1\}, \quad \forall i \in I, j \in J, m \in M(i), \label{eq:full_binary_cons3} \\
& y_{im} \in \{ 0, 1\}, \quad \forall i \in I, m \in M(i), \label{eq:full_binary_cons2} \\
& x_{ij} \in \{ 0, 1\}, \quad \forall i \in I, j \in J. \label{eq:full_binary_cons}
\end{align}
\end{subequations} 
This formulation follows the same structure as the DAO model presented in Section \ref{sec_DAO}. The constraints of $\mathcal{X}$ are now explicitly defined by constraint set \eqref{eq:full_x_1}, which ensures that all tasks are assigned to exactly one underwriter, and constraint set \eqref{eq:full_x_2}, which enforces that the number of tasks from $J_t$ assigned to each worker does not exceed their daily capacity, $c_{it}$, for any $t \in T$. As noted previously, when we only consider one candidate predictive model per worker ($M(i) = 1$), the formulation simplifies to the UAO model presented in Section \ref{sec_UAO}.

\section{Case Study Results} \label{sec_results}
In this section, we present the findings of the case study. First, in Section~\ref{cs_obj_val}, we examine the expected value of the solutions generated by the three UAO approaches and the proposed DAO approach when optimizing underwriter allocation. This analysis focuses on how each method for leveraging predictive insights impacts the profitability of the allocation strategy. Next, in Section~\ref{sec_vary_data}, we investigate how varying the amount of training data influences optimization outcomes, highlighting the adaptability of each method across different data settings and the data-responsiveness of the DAO framework.  In Section~\ref{sec_dao_selection}, we examine the predictive model selection of the DAO approach, exploring its implications for workforce optimization and offering insights that could inform other applications. Finally, in Section~\ref{sec_worker_allocation}, we analyze how tasks are allocated within the DAO solution, focusing on the distribution of task values and the strategic use of different predictive models based on task characteristics.

\subsection{Worker Allocation Objective Value}\label{cs_obj_val}
\begin{figure}[ht!]
    \centering
    \includegraphics[width = \textwidth]{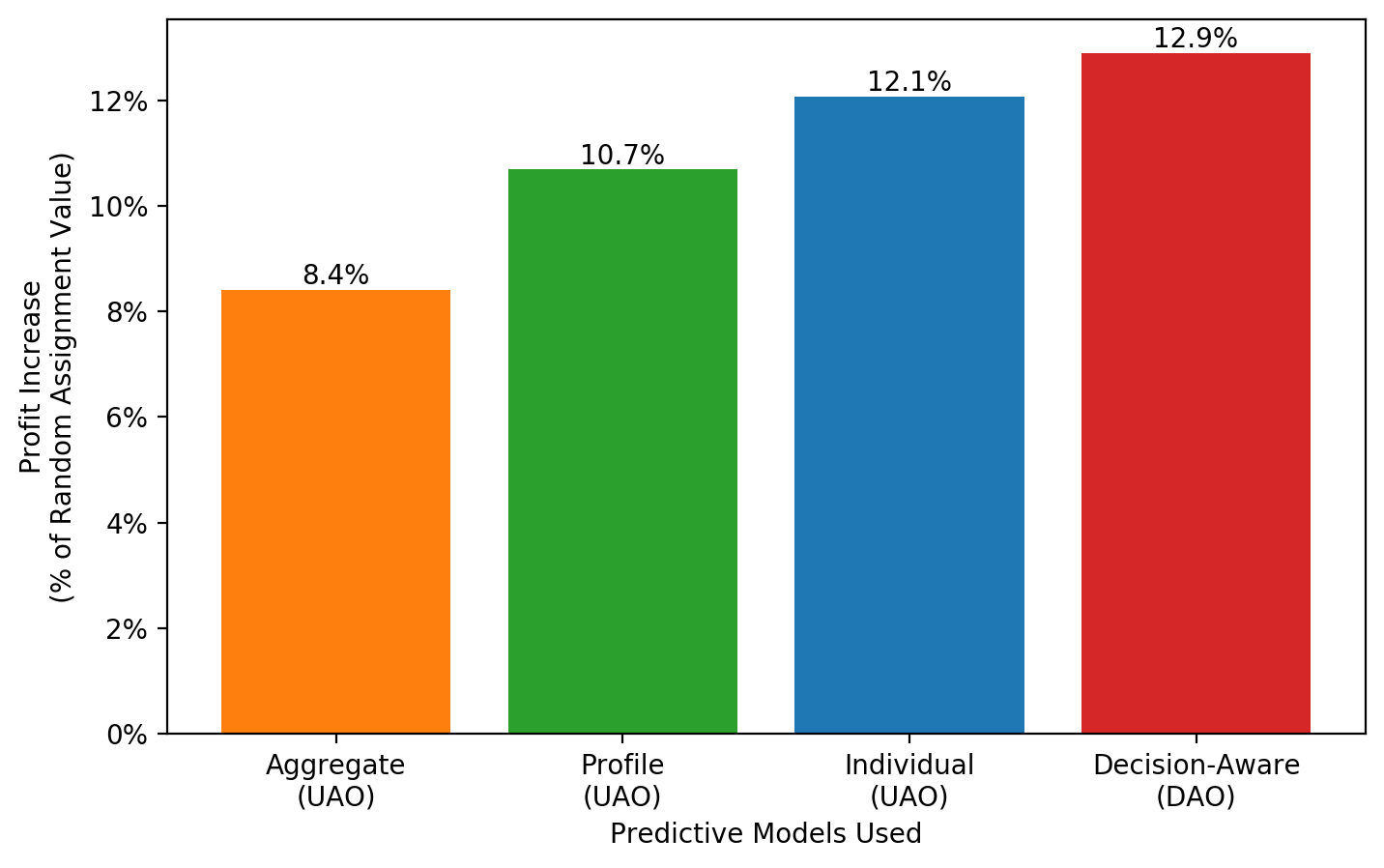}
    \caption{Expected profit increase compared to random worker allocation when different behavioral predictive models are used in the optimization framework.} \label{fig:exp_value_comp}
\end{figure}

In our case study, we evaluated three distinct UAO approaches, each leveraging a different predictive model to forecast underwriter behavior, as well as the DAO approach, where the predictive model representing each underwriter’s behavior is dynamically selected within the allocation model. The goal was to assess how each approach, rooted in different predictive strategies, impacts the profitability of workforce allocation. The expected value of the resulting worker allocations, computed as the-percentage improvement over the expected value of the random policy, is presented in Figure \ref{fig:exp_value_comp}. The expected value of a random policy is computed as the expected predicted value obtained by randomly assigning tasks to workers and predictive models, subject to the allocation constraints.

The results indicate that all four optimization approaches significantly outperformed the randomized assignment policy in terms of expected profitability. Among the UAO approaches, which rely on a single predictive model for each underwriter, the aggregate model led to an 8.4-percent increase in profits, the profile model yielded a 10.7-percent increase, and the individualized model achieved a 12.1-percent increase. Interestingly, despite the aggregate model having the highest AUC (Figure~\ref{fig:ROC}), it produced the smallest profit improvement, while the individualized model, with the lowest AUC, delivered the highest profit increase. This highlights a key insight: while overall predictive accuracy (as indicated by AUC) is important, the ability to capture individual nuances in worker behavior is crucial for optimizing workforce allocation.

The DAO approach, which allows for the selection of the optimal predictive model for each underwriter within the optimization process, demonstrated the best overall performance, achieving a 12.9-percent improvement over the randomized policy. This prompts a key consideration as to whether the additional 0.8-percent profit gain from the DAO approach is worth the added complexity of implementing a more sophisticated model selection strategy. In many cases, the individualized models within the UAO framework captured the majority of the potential performance improvements, suggesting that this approach may be more suitable for organizations that prioritize operational simplicity over the complexities of the DAO methodology.

\subsection{Impact of Varying Training Data}\label{sec_vary_data}
To further explore the adaptability of each approach, we analyzed how varying the amount of training data influences the performance of the UAO and DAO approaches. This analysis assesses how both the quantity and recency of data impacts predictive accuracy and the resulting optimization outcomes. This is an important consideration in industries like auto insurance, where conditions can change rapidly.

\begin{figure}[ht!]
    \centering
    \includegraphics[width = \textwidth]{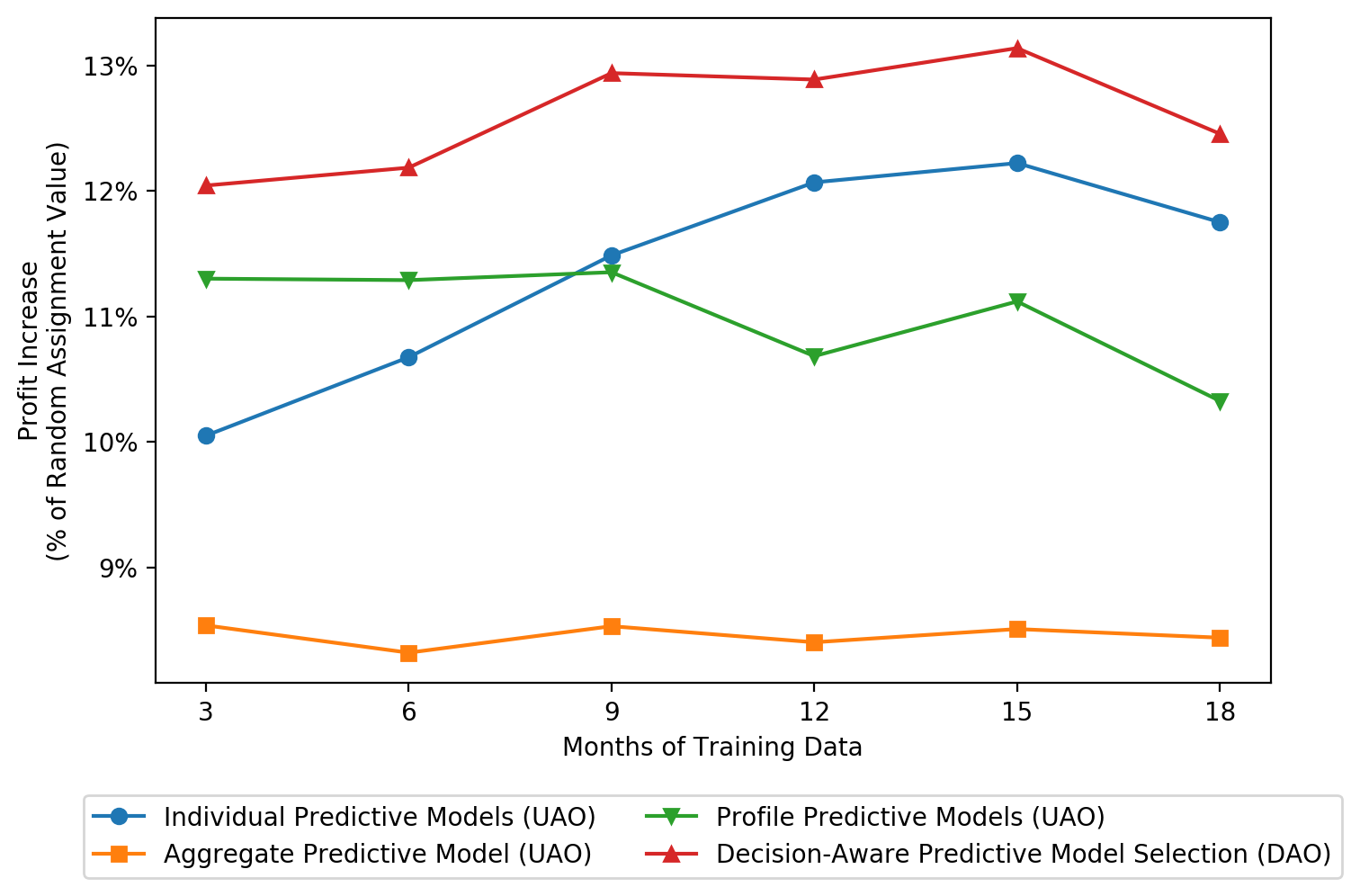}
    \caption{Expected profit increase compared to random worker allocation for different predictive model selections and varying amounts of training data.} \label{fig:exp_value_varying_data}
\end{figure}

We evaluated the models using training data ranging from the most recent 3 to 18 months, in 3-month increments to understand the effects of data recency. Figure \ref{fig:exp_value_varying_data} illustrates the impact of these different training data sets on the expected value of the resulting worker allocations, computed as the-percentage improvement over the expected value of the random policy.

The results reveal that the UAO approach using the aggregate predictive model consistently delivered an 8-percent profit increase, regardless of data size. The UAO approach with profile predictive models performed better, achieving approximately an 11-percent profit increase. With smaller datasets (3-6 months), the UAO profile model proved to be the most effective, but as more training data became available, the individualized model emerged as the top performer, reaching a 12-percent profit increase with 15 months of training data.

In all scenarios, the DAO approach outperformed the UAO strategies, with its largest improvement resulting in a 1.5-percent additional profit increase over the best-performing UAO strategy. This occurred when predictive models were trained on 9 months of data. We hypothesize that the DAO approach offers the greatest improvements over UAO at this point because, around the 9-month range, the benefit of increasing data on AUC began to plateau for all three predictive models (Appendix C). This suggests that the DAO approach is most effective when the predictive models have all reached a reasonable level of predictive power. This is likely because the models need to be sufficiently reliable to make the dynamic selection process meaningful, enabling the DAO framework to leverage the strengths of each predictive model and optimize workforce allocation effectively.

These findings align with the underlying intuition of the problem. Although the aggregate predictive model excelled in general predictive accuracy (Figure \ref{fig:ROC}), it struggled to capture the specific nuances of individual worker behavior necessary for optimal assignment policies. The profile model was particularly effective with limited data, balancing general worker tendencies with sufficient data for accuracy. As the training data increased, the individualized models gained a distinct advantage, overcoming the limitations of smaller datasets. Ultimately, the DAO approach demonstrated its strength by adapting flexibly and efficiently to a range of data conditions, making it the most robust strategy across both data-scarce and data-rich environments.

\subsection{DAO Predictive Model Selection} \label{sec_dao_selection}
\begin{figure}[ht!]
    \centering
    \includegraphics[width = .8\textwidth]{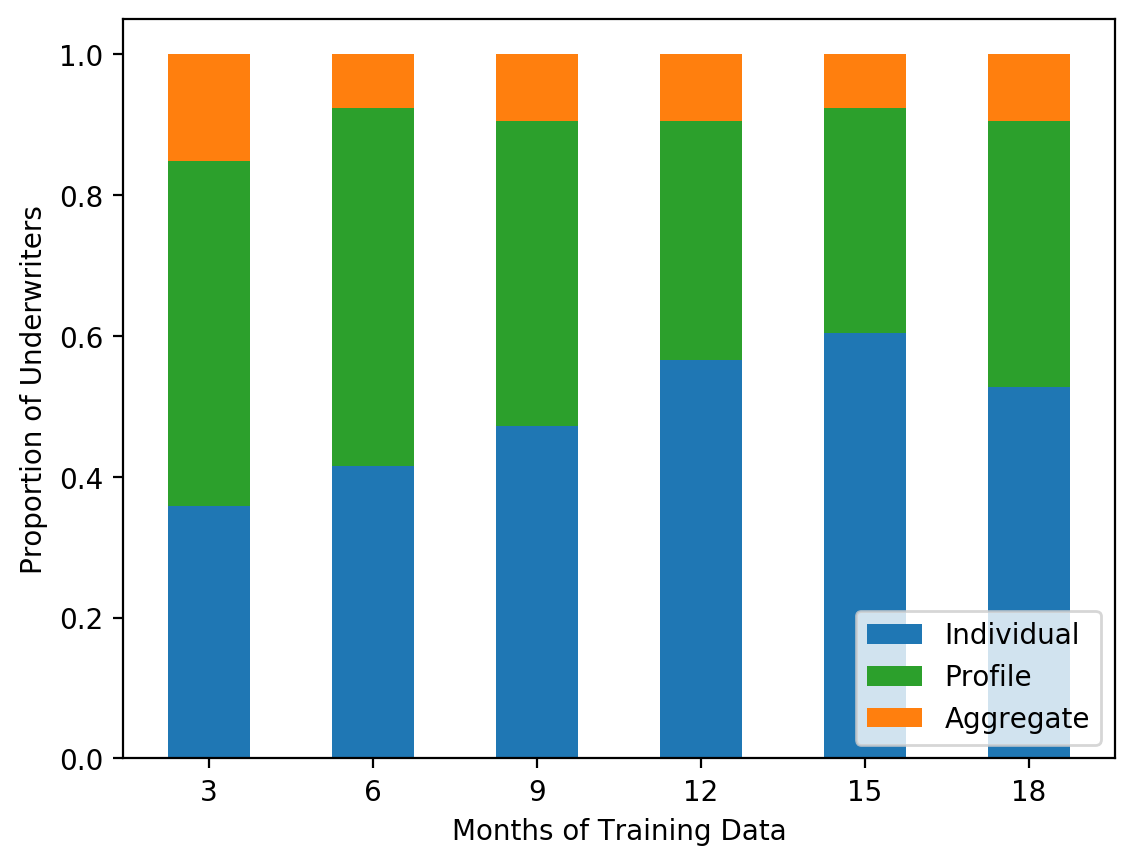}
    \caption{Predictive model selection in the DAO approach across varying amounts of training data.} \label{fig:predictive_model_selection}
\end{figure}

In this section, we examine the insights provided by the DAO predictive model selection. Figure \ref{fig:predictive_model_selection} shows the predictive model choices within the DAO approach across varying amounts of training data.

As shown, with limited data (3 to 6 months), the profile model is predominantly selected for most workers. However, as more data becomes available, the preference shifts towards individualized predictive models. This shift makes sense because individual models require more data to yield accurate results. The key takeaway from this analysis is that the aggregate model is rarely the best choice for predicting worker behavior. Even when predictive power is limited, it is better to leverage the nuances of worker behavior that are captured through more personalized predictive models. When data is scarce, the profile model effectively represents workers, but as more data accumulates, transitioning to individualized models provides a more precise understanding of each worker's behavior, which can be better utilized in assignment policies.

\subsection{Analysis of Worker Allocation Within the DAO Solution} \label{sec_worker_allocation}

In this section, we examine how the predictive models influence the allocation of tasks to workers. This analysis sheds light on the decision-making process within the DAO solution, particularly in how task values are distributed across workers and the role of different predictive models in these assignments. In our analysis, we examine how workers are allocated to tasks within the DAO solution, focusing on a scenario with 12 months of training data. We first analyze the value of the tasks assigned to each worker. 

\begin{figure}[ht!]
    \centering
    \includegraphics[width = \textwidth]{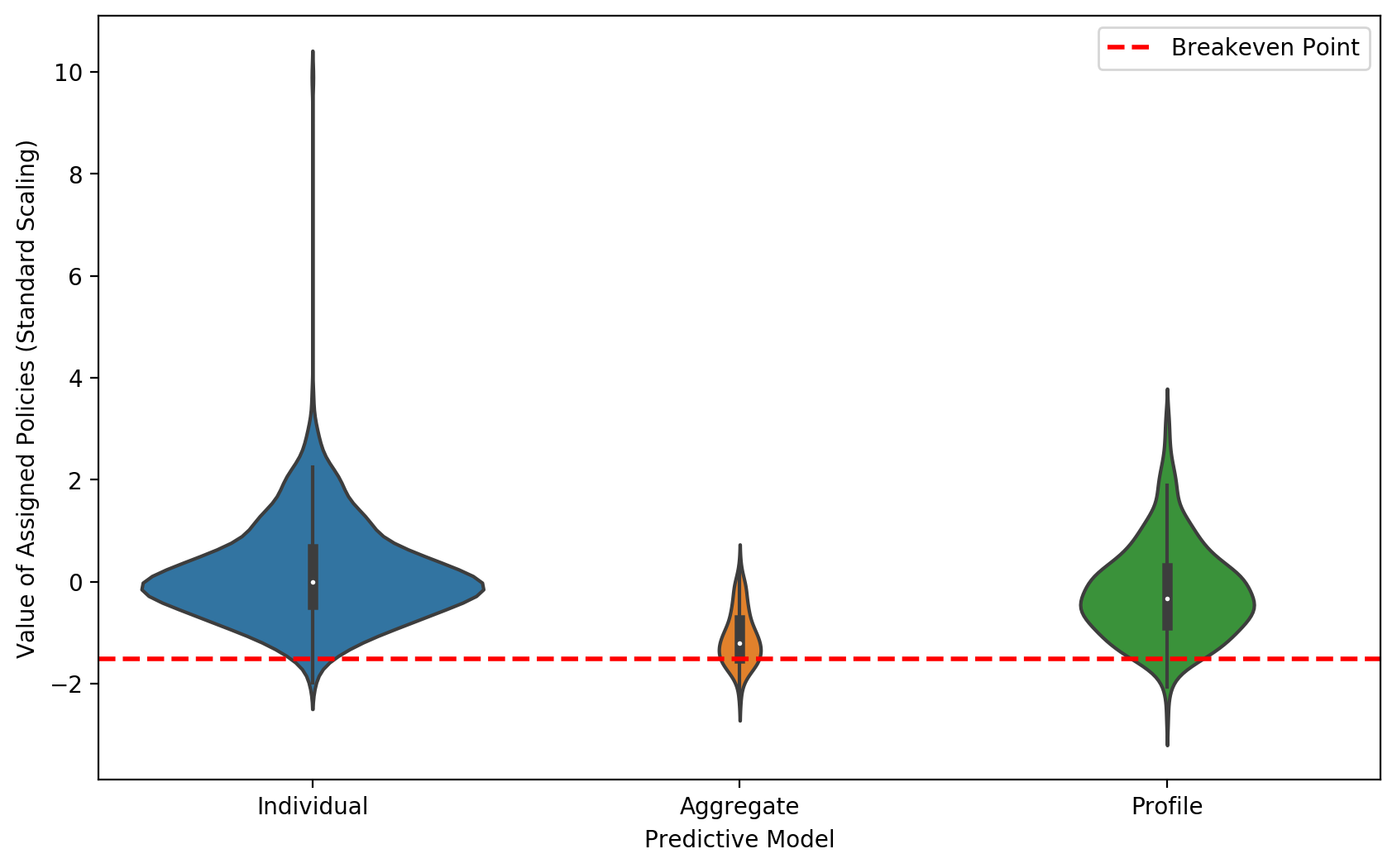}
    \caption{Distribution of task values assigned to workers using individual, profile, and aggregate predictive models within the DAO framework. The red dashed line indicates policies where \(v_j = 0\).} \label{fig:violin}
\end{figure}

Figure~\ref{fig:violin} displays violin plots that show the value of tasks that are assigned to workers using either the individual, aggregate, or profile model. Note that the individual model violin plot is much larger than the aggregate model because more tasks are assigned using the individual model (see Figure~\ref{fig:predictive_model_selection}). The median value (after standard scaling) for the individual, aggregate, and profile policies is 0.001, -1.208, and -0.328, respectively. The distribution of assigned tasks shows that higher-value tasks are typically allocated to workers represented by individualized models. Tasks assigned using profile models tend to have slightly lower values, while tasks assigned based on the aggregate model often have values closer to or even below the breakeven point. This observation highlights a clear difference in the value of tasks allocated depending on the predictive model used within the DAO framework.

\begin{figure}[ht!]
    \centering
    \includegraphics[width = \textwidth]{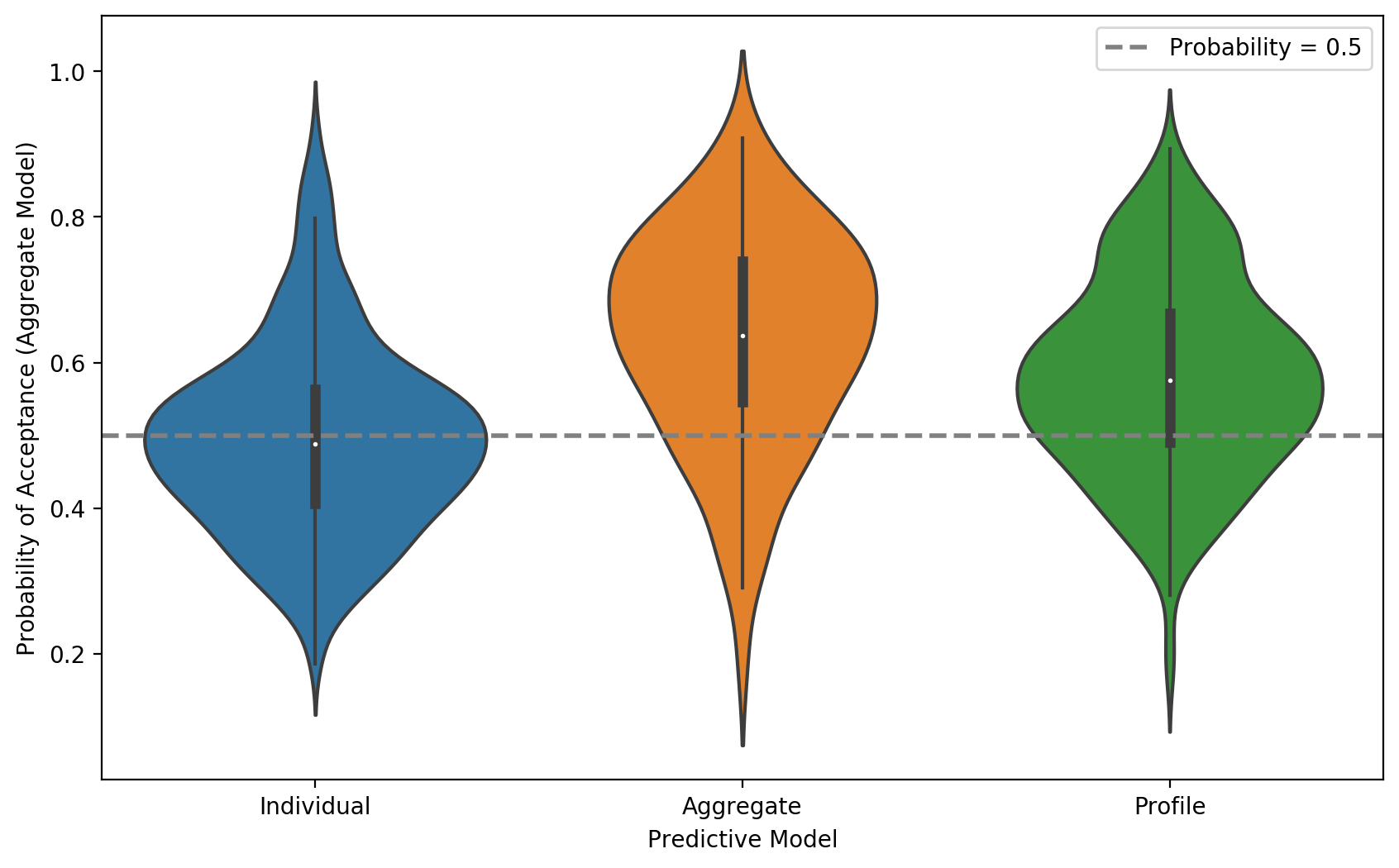}
    \caption{Aggregate model predictions for tasks assigned to the workers represented by individual, aggregate, or profile models.} \label{fig:violin_var}
\end{figure}

To further understand the model assignments, we take a closer look at the acceptance probabilities of the assigned tasks. For each task, we obtain the aggregate model prediction for that task (note this value is the same for all workers). Figure~\ref{fig:violin_var} displays violin plots of the aggregate model prediction partitioned by tasks that were optimally assigned using either the individual, aggregate, or profile models.  The median predicted probability made by the aggregate predictive model for tasks optimally assigned using the individual, aggregate, and profile models is 0.489, 0.637, and 0.576, respectively. 

A predicted probability near 0.5 indicates that aggregate predictive model has low confidence in it's ability to predict worker behavior for that task. The first violin shows that when the aggregate model predictions are close to 0.5, the optimization model tends to choose individual model predictions. A similar pattern is observed for the profile violin plot. Conversely, the second violin shows that the aggregate model is most often selected in the optimization when its predicted probabilities are high, indicating aggregate model confidence and suggesting consensus across workers. This pattern highlights two key insights. First, individual and profile models are often employed when there are greater individual differences in the predictive models. This approach makes sense, as it allows the predictive model selection to effectively leverage these differences. Second, tasks with higher variance likely have higher values, which aligns with the nature of the insurance business, where policies with less consensus among workers tend to involve higher risk and may potentially lead to greater profitability (Appendix D).

These results suggest that the DAO framework strategically assigns higher-risk or more uncertain tasks to workers using individualized models, which aim to capture subtle differences in worker behavior. Conversely, tasks where worker responses are more predictable are assigned using aggregate or profile models. This indicates that employing a more personalized, individualized predictive modeling approach may be beneficial when greater variance in worker behavior is expected, even if predictive performance metrics like AUC are lower than those of less individualized models trained with more data.

\section{Conclusion} \label{sec_conclusion}
This paper introduces a novel approach to workforce allocation that leverages predictive modeling and optimization to enhance worker-task assignments. By dynamically selecting the optimal predictive model for each worker within the proposed DAO framework, we improve alignment between individual worker characteristics, predictive modeling techniques, and organizational objectives. The DAO framework not only enhances decision-making accuracy but also adapts flexibly to various data environments, effectively meeting diverse organizational needs.

We present a case study based on real-world data from an auto-insurance company, offering key insights into optimizing workforce management through the application of various predictive modeling approaches within UAO and DAO frameworks. By evaluating individual, aggregate, and profile models, our analysis demonstrates that predictive insights can significantly improve workforce allocation. The DAO framework, in particular, stands out for its ability to dynamically select the most appropriate model for each worker, aligning the optimization process with specific organizational goals. This flexibility makes the DAO approach especially advantageous in scenarios with fluctuating data conditions or when the optimal predictive model is uncertain or dependent on specific worker assignments. Our findings highlight the potential of predictive modeling to refine decision-making in workforce management, driving better outcomes in complex, real-world environments.

For future research, several promising avenues exist to advance the integration of predictive models within optimization frameworks. One potential direction is to extend our proposed DAO framework by replacing the expected value objective function with stochastic programming formulations, such as Conditional Value-at-Risk (CVaR). This could provide interesting insights into how predictive model selection changes in different risk environments and with other sources of uncertainty. Additionally, we encourage future work that explores decision-aware strategies in others stages of the machine learning to optimization pipeline. By critically assessing the integration of predictive models within optimization frameworks, researchers and practitioners can develop more resilient, effective, and adaptable strategies to better address the complexities of real-world scenarios.

\ACKNOWLEDGMENT{Support for this research was provided by American Family Insurance through a research partnership with the University of Wisconsin–Madison’s Data Science Institute.}

\bibliographystyle{informs2014} 
\bibliography{sample}

\end{document}


\section*{Supplementary Material to ``Decision-Aware Predictive Model Selection for Workforce Allocation''}
\text{ }
 \begin{APPENDICES}
        \numberwithin{figure}{section}
        \numberwithin{table}{section}
        \renewcommand{\thefigure}{\thesection.\arabic{figure}}
        \renewcommand{\thetable}{\thesection.\arabic{table}}
    \section{Data Inclusion Diagram}\label{sec:data_inclusion}
    The data inclusion process began with a dataset of 48,284 auto-insurance policies reviewed by underwriters during the study period. Policies reviewed by multiple underwriters---common practice when training junior underwriters---were excluded, resulting in 47,349 policies reviewed by a single underwriter. To ensure an accurate assessment of each predictive model's fit, we further excluded policies reviewed by underwriters who had evaluated fewer than 50 policies in the test dataset, leaving 37,782 policies. These policies were then divided into a training set, a testing set (3 months, $n = 3,857$), and an optimization set (1 month, $n = 1,187$). The size of the training dataset varied depending on the number of months used for model training (as discussed in Section 5.2), ranging from 3 months with 3,752 policies to 18 months with 32,738 policies. The data inclusion process is summarized in Figure~\ref{fig:data-flow-chart}. The amount of data used to train and test each predictive model varied based on the underwriter the model was intended to represent and the specific data strategy of the predictive modeling approach (individual, aggregate, or profile).

    \begin{figure}[ht]
    \centering
    \begin{tikzpicture}[font=\footnotesize] 
    \draw[thick] (0,0) rectangle ++(3.5,1.75) node[pos=.5, align=center] {Training Data Set \\ 3-18 months \\ ($n = 3,752 - 32,738$)};
    \draw[thick] (4,0) rectangle ++(3.5,1.75) node[pos=.5, align=center] {Testing Data Set \\ 3 months \\ ($n = 3,857$)};
    \draw[thick] (8,0) rectangle ++(3.5,1.75) node[pos=.5, align=center] {Optimization Data Set \\ 1 month \\ ($n = 1,187$)};
    
    \draw[thick] (3,2.5) rectangle ++(5.5,1.75) node[pos=.5, align=center] {Reviewing underwriter assigned \\ 50+ policies in testing data set.\\ ($n = 37,782$)};
    \draw[thick] (3,5) rectangle ++(5.5,1.75) node[pos=.5, align=center] {Reviewed by one underwriter. \\ ($n = 47,349$)};
    \draw[thick] (3,7.5) rectangle ++(5.5,1.75) node[pos=.5, align=center] {Auto-insurance policies reviewed by\\an underwriter during the study period. \\ ($n = 48,284$)};
    
    \draw[thick, -latex] (5.75,2.5) -- (5.75,1.75);
    \draw[thick, -latex] (5.75,5) -- (5.75,4.25);
    \draw[thick, -latex] (5.75,7.5) -- (5.75,6.75);
    
    \draw[thick, -latex] (1.75,2.125) -- (1.75,1.75);
    \draw[thick] (1.75,2.125) -- (5.75,2.125);
    \draw[thick, -latex] (9.75,2.125) -- (9.75,1.75);
    \draw[thick] (9.75,2.125) -- (5.75,2.125);

    \end{tikzpicture}
    \caption{Data inclusion diagram for auto-insurance case study.}
    \label{fig:data-flow-chart}
    \end{figure}
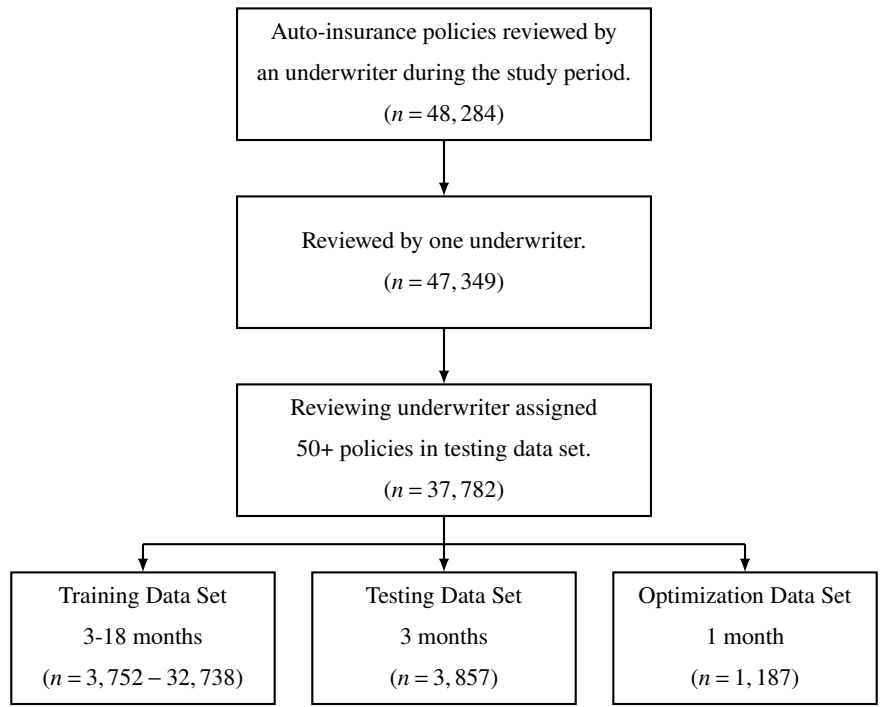
 
   \section{Profile Predictive Model Clustering} \label{sec:app_cluster}
    In the profile predictive model approach, we grouped workers into clusters. When working with datasets spanning six months or longer, we identified six distinct clusters. However, for experiments with only three months of data, we opted for four clusters. This clustering strategy was informed by analyzing the distribution of workers relative to the aggregate model, coupled with a detailed examination of intercluster distances. These distances are visualized in Figure \ref{fig:elbow}, with the x-axis indicating the number of clusters considered and the y-axis showing inertia—a metric representing the sum of squared distances between each data point and its respective cluster's centroid. This metric helped gauge the cohesiveness of the clusters. An observable elbow in the graph around six clusters guided the decision for longer datasets, while four clusters were deemed more appropriate for the shorter, 3-month data scenario, achieving similar inertia scores and facilitating more coherent groupings. This approach ensures the predictive models are grounded in a nuanced understanding of worker behavior patterns over varying time frames.

    \begin{figure}[ht!]
    \centering
    \includegraphics[width = 0.80\textwidth]{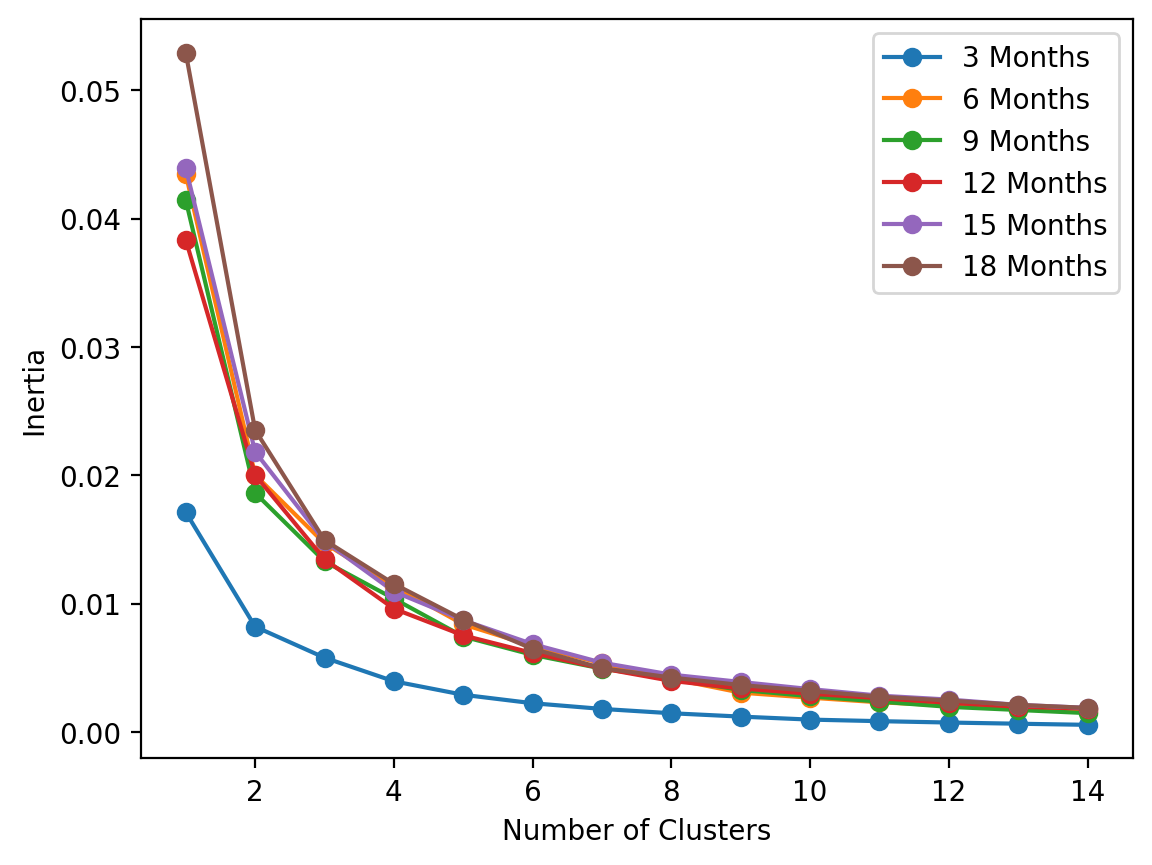}
    \caption{Inertia versus cluster numbers for profile predictive models for various amounts of training data.} \label{fig:elbow}
    \end{figure}

    \section{Predictive Model Accuracy} \label{sec:app_acc}
    Figure \ref{fig:auc} details the test set AUC of the predictive models in relation to the volume of training data used. The graph's x-axis represents the duration of the training data in months, while the y-axis measures the weighted (by number of observations) average AUC of the predictive models. We calculated AUC for individual and profile predictive models using the predictions corresponding to the worker who reviewed each policy. As shown, all three plots exhibit an elbow around the 9-month mark, indicating that the benefit of adding more data starts to plateau beyond this point. This suggests that while increasing the volume of training data initially improves model performance, further increases yield diminishing returns. The plateau may be attributed to the models having already captured the most relevant patterns in the data, after which additional data contributes less to improving predictive power.

    \begin{figure}[ht!]
        \centering
        \includegraphics[width = \textwidth]{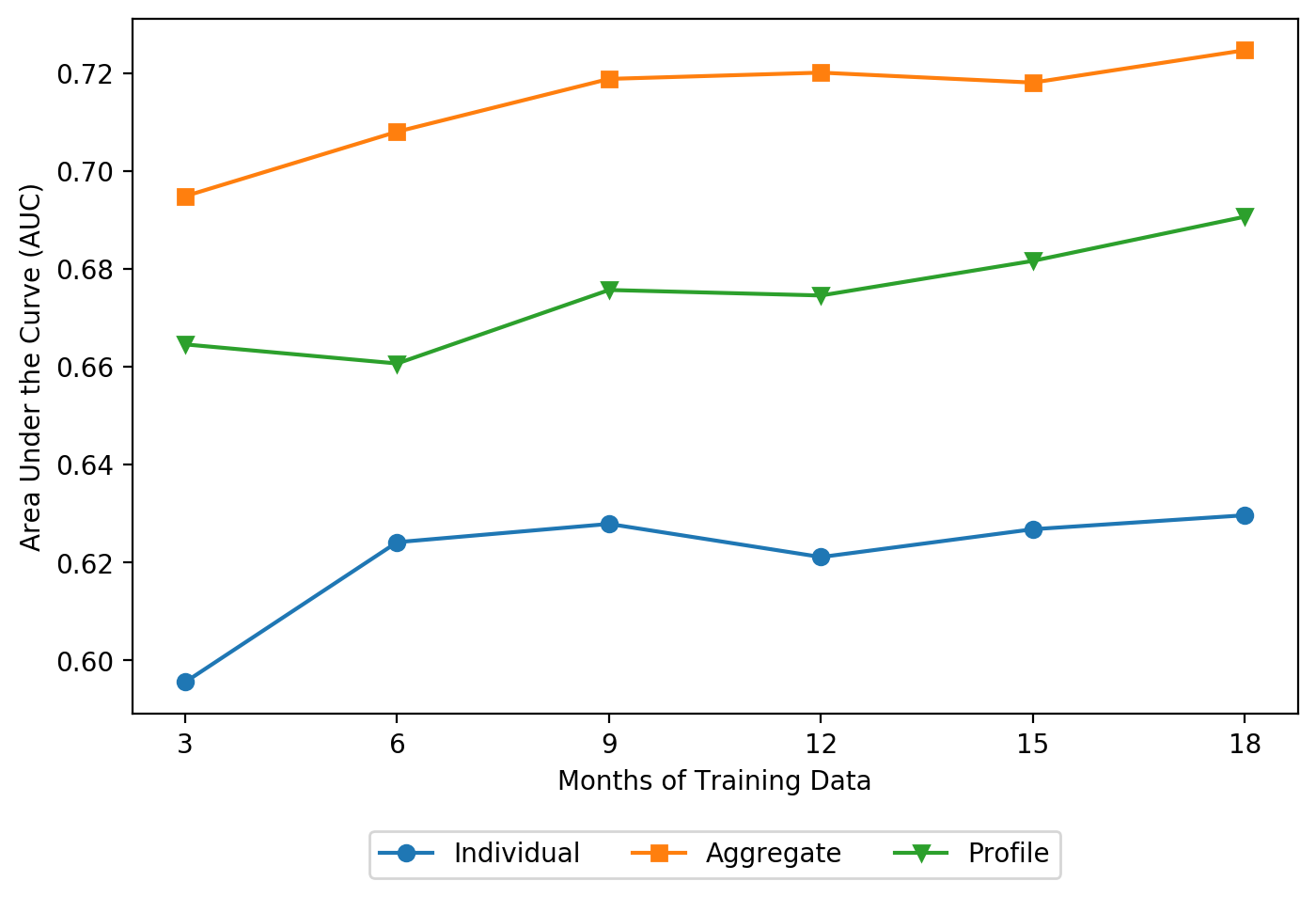}
        \caption{Average AUC of the predictive models for various amounts of training data.} \label{fig:auc}
    \end{figure}

\section{Consensus and Policy Value Analysis} \label{sec:appendix_consensus_value}
Figure~\ref{fig:value_variance} shows the relationship between the variance in predicted worker behavior from the individual predictive models and the policy value. Each point represents a policy. The x-axis displays the variance in predicted worker behavior, reflecting how much the predictions differ across the individual predictive models, with higher values indicating less consensus. The y-axis represents the policy value, with higher values indicating policies that are more valuable relative to the average, based on the standardized scale.

As shown, higher variance among the predictive models tends to correspond with higher policy values. This suggests that policies with less consensus—indicative of greater uncertainty—may be associated with higher profitability. This aligns with the nature of the insurance business, where higher-risk policies, as indicated by greater variance in predictions, tend to correspond with higher policy values.

\begin{figure}[ht!]
    \centering
    \includegraphics[width=00.80\textwidth]{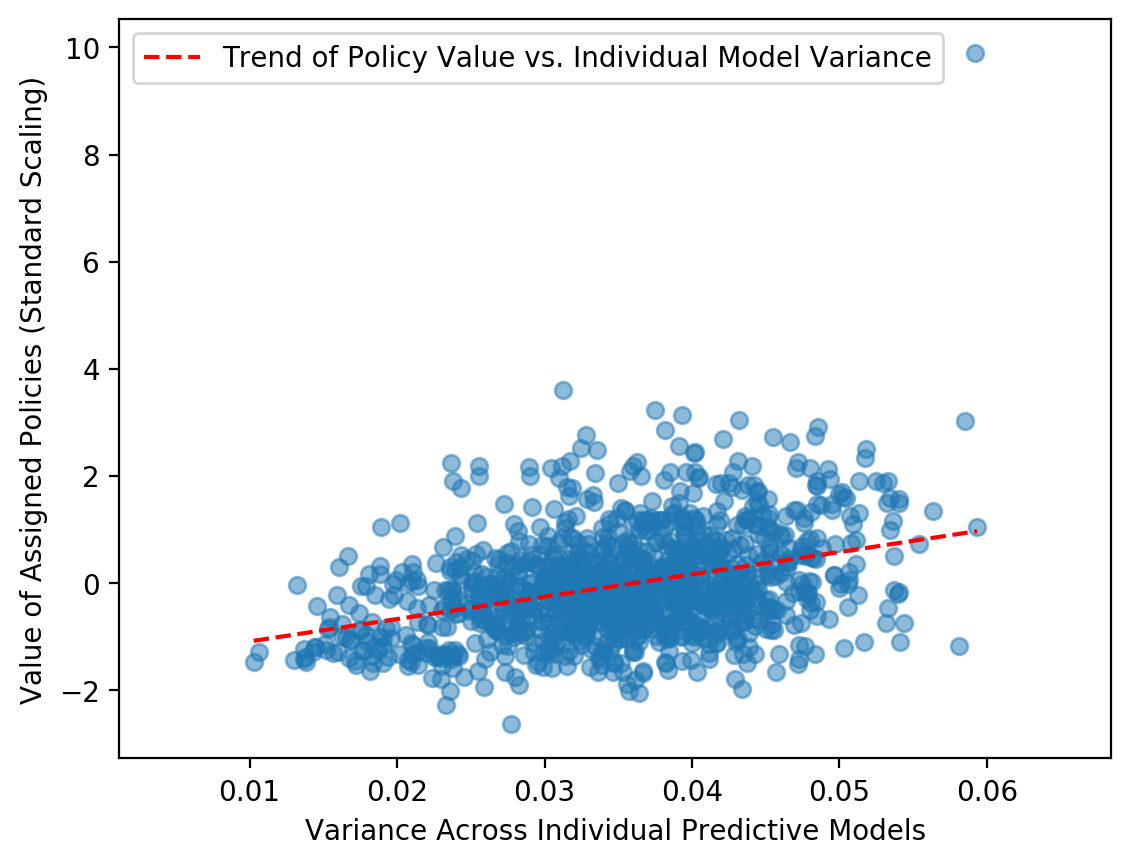}
    \caption{Relationship between the variance in predicted worker behavior and policy value. Each point represents a policy. The dashed red line indicates the observed trend, suggesting that policies with higher variance in individual worker behavior tend to correspond with higher policy values.}
    \label{fig:value_variance}
\end{figure}

 \end{APPENDICES}